\documentclass[english,12pt]{smfart}
\usepackage{amssymb}
\usepackage{amsfonts}
\usepackage{amscd}
\usepackage[all]{xypic}

\CompileMatrices

\swapnumbers

\def\ord{{\operatorname {\rm ord}}}
\def\ac{{\overline{\operatorname  {\rm ac}}}}

\def\Def{{\rm Def}}

\def\RDef{{\rm RDef}}
\def\RDefe{{\rm RDef}^{\rm exp}}
\def\RDefw{{\rm RDef}^{\rm e}}

\def\LPre{\cL_{\rm DP,P}}

\let\cal\mathcal

\def\11{{\mathbf 1}}
\def\ordjac{{\operatorname {\rm ordjac}}}
\def\llp{\mathopen{(\!(}}
\def\llb{\mathopen{[\![}}
\def\rrp{\mathopen{)\!)}}
\def\rrb{\mathopen{]\!]}}

\mathchardef\alphag="7C0B
\mathchardef\betag="7C0C
\mathchardef\gammag="7C0D
\mathchardef\deltag="7C0E
\mathchardef\varepsilong="7C22
\mathchardef\varphig="7C27
\mathchardef\psig="7C20
\mathchardef\zetag="7C10
\mathchardef\epsilong="7C0F
\mathchardef\rhog="7C1A
\mathchardef\taug="7C1C
\mathchardef\upsilong="7C1D
\mathchardef\iotag="7C13
\mathchardef\thetag="7C12
\mathchardef\pig="7C19
\mathchardef\sigmag="7C1B
\mathchardef\etag="7C11
\mathchardef\omegag="7C21
\mathchardef\kappag="7C14
\mathchardef\lambdag="7C15
\mathchardef\mug="7C16
\mathchardef\xig="7C18
\mathchardef\chig="7C1F
\mathchardef\nug="7C17
\mathchardef\varthetag="7C23
\mathchardef\varpig="7C24
\mathchardef\varrhog="7C25
\mathchardef\varsigmag="7C26
\mathchardef\Omegag="7C0A
\mathchardef\Thetag="7C02
\mathchardef\Sigmag="7C06
\mathchardef\Deltag="7C01
\mathchardef\Phig="7C08
\mathchardef\Gammag="7C00
\mathchardef\Psig="7C09
\mathchardef\Lambdag="7C03
\mathchardef\Xig="7C04
\mathchardef\Pig="7C05
\mathchardef\Upsilong="7C07

\newtheorem{theorem}[subsubsection]{Theorem}
\newtheorem{lem}[subsubsection]{Lemma}

\newtheorem{prop}[subsubsection]{Proposition}

\newtheorem{claim}[subsubsection]{Claim}

\theoremstyle{definition}
\newtheorem{definition}[subsubsection]{Definition}

\newtheorem{def-prop}[subsubsection]{Proposition-Definition}
\newtheorem{def-theorem}[subsubsection]{Theorem-Definition}
\newtheorem{def-lem}[subsubsection]{Lemma-Definition}

\theoremstyle{remark}
\newtheorem{remark}[subsubsection]{Remark}

\theoremstyle{plain}

\numberwithin{equation}{subsection}

\def\boxit#1#2{\setbox1=\hbox{\kern#1{#2}\kern#1}%
\dimen1=\ht1 \advance\dimen1 by #1
\dimen2=\dp1 \advance\dimen2 by #1
\setbox1=\hbox{\vrule height\dimen1 depth\dimen2\box1\vrule}%
\setbox1=\vbox{\hrule\box1\hrule}%
\advance\dimen1 by .4pt \ht1=\dimen1
\advance\dimen2 by .4pt \dp1=\dimen2 \box1\relax}

\newcommand{\sur}[2]{\genfrac{}{}{0pt}{}{#1}{#2}}
\renewcommand{\theequation}{\thesubsection.\arabic{equation}}

\let\cal\mathcal

\def\AA{{\mathbf A}}

\def\CC{{\mathbf C}}

\def\ee{{\mathbf e}}
\def\FF{{\mathbf F}}
\def\GG{{\mathbf G}}

\def\LL{{\mathbf L}}

\def\NN{{\mathbf N}}

\def\QQ{{\mathbf Q}}

\def\ZZ{{\mathbf Z}}

\def\cA{{\mathcal A}}
\def\cB{{\mathcal B}}
\def\cC{{\mathcal C}}
\def\cD{{\mathcal D}}

\def\cF{{\mathcal F}}

\def\cI{{\mathcal I}}

\def\cL{{\mathcal L}}

\def\cO{{\mathcal O}}
\def\cP{{\mathcal P}}

\def\cS{{\mathcal S}}

\def\cU{{\mathcal U}}

\mathchardef\alphag="7C0B
\mathchardef\betag="7C0C
\mathchardef\gammag="7C0D
\mathchardef\deltag="7C0E
\mathchardef\varepsilong="7C22
\mathchardef\varphig="7C27
\mathchardef\psig="7C20
\mathchardef\zetag="7C10
\mathchardef\epsilong="7C0F
\mathchardef\rhog="7C1A
\mathchardef\taug="7C1C
\mathchardef\upsilong="7C1D
\mathchardef\iotag="7C13
\mathchardef\thetag="7C12
\mathchardef\pig="7C19
\mathchardef\sigmag="7C1B
\mathchardef\etag="7C11
\mathchardef\omegag="7C21
\mathchardef\kappag="7C14
\mathchardef\lambdag="7C15
\mathchardef\mug="7C16
\mathchardef\xig="7C18
\mathchardef\chig="7C1F
\mathchardef\nug="7C17
\mathchardef\varthetag="7C23
\mathchardef\varpig="7C24
\mathchardef\varrhog="7C25
\mathchardef\varsigmag="7C26
\mathchardef\Omegag="7C0A
\mathchardef\Thetag="7C02
\mathchardef\Sigmag="7C06
\mathchardef\Deltag="7C01
\mathchardef\Phig="7C08
\mathchardef\Gammag="7C00
\mathchardef\Psig="7C09
\mathchardef\Lambdag="7C03
\mathchardef\Xig="7C04
\mathchardef\Pig="7C05
\mathchardef\Upsilong="7C07

\DeclareMathOperator*{\sq}{\square}
\def\ord{{\operatorname {\rm ord}}}
\def\Jac{{\operatorname  {\rm Jac}}}
\def\Fv{{\mathfrak{F}}}
\def\Fr{{\mathfrak{f}}}

\def\11{{\mathbf 1}}
\def\ordjac{{\operatorname {\rm ord jac}}}

\def\LO{{\cL_{\cO}}}

\def\Def{{\rm Def}}

\def\Def{{\rm Def}}

\def\RDef{{\rm RDef}}

\let\cal\mathcal

\def\LPre{\cL_{\rm DP,P}}

\begin{document}

\title[Constructible exponential functions]{Constructible exponential functions,
motivic Fourier transform and transfer principle}

\author{Raf Cluckers}
\address{Katholieke Universiteit Leuven,
Departement wiskunde, Celestijnenlaan 200B, B-3001 Leu\-ven,
Bel\-gium. Current address: \'Ecole Normale Sup\'erieure,
D\'epartement de ma\-th\'e\-ma\-ti\-ques et applications, 45 rue
d'Ulm, 75230 Paris Cedex 05, France} \email{cluckers@ens.fr}
\urladdr{www.dma.ens.fr/$\sim$cluckers/}

\author{Fran\c cois Loeser}

\address{{\'E}cole Normale Sup{\'e}rieure,
D{\'e}partement de math{\'e}matiques et applications, 45 rue
d'Ulm, 75230 Paris Cedex 05, France (UMR 8553 du CNRS)}
\email{Francois.Loeser@ens.fr} \urladdr{www.dma.ens.fr/~loeser/}


\maketitle

\renewcommand{\partname}{}

\def\cC{{\cal C}}

\section{Introduction}
In our previous
work \cite{cl}, 
we laid general foundations
for motivic integration of constructible functions.
One of the most salient features of  motivic
constructible functions is that they form a class which is stable by direct image and
that motivic integrals of
constructible functions depending on parameters are constructible as functions of the parameters. Though motivic constructible functions as defined
in  \cite{cl} encompass motivic analogues of many  functions occuring
in integrals over
non archimedean local fields, one important class of functions was still missing in the picture,
namely motivic analogues  of
non archimedean integrals of the type
\begin{equation*}
\int_{\QQ_p^n} f (x) \Psi (g (x)) \vert dx\vert,
\end{equation*}
with $ \Psi$ a (non trivial) additive character
on $\QQ_p$, $f$ a $p$-adic constructible function and $g$ a
$\QQ_p$-valued definable function on $\QQ_p^n$, and their
parametrized versions, functions of the type
\begin{equation*}
\lambda \longmapsto
\int_{\QQ_p^n} f (x, \lambda) \Psi (g (x, \lambda)) \vert dx\vert,
\end{equation*}
where $\lambda$ runs over, say $\QQ_p^m$, and
$f$ and $g$ are now functions on $\QQ_p^{m + n}$.
Needless to say, such kind of integrals are ubiquitous in harmonic analysis
over
non archimedean local fields,
$p$-adic representation Theory and
the Langlands Program.

\medskip

One of the  purposes of the present paper is to fill this gap by
extending the framework of \cite{cl} in order to include motivic
analogues of exponential integrals of the above type. Once this is
done one is able to develop a natural Fourier transform and to prove
various forms of Fourier inversion. Another interesting feature of
our  formalism is that it makes possible to state and prove a
general transfer principle for integrals over non archimedean local
fields,  allowing to transfer identities between  functions defined
by integrals over fields of characteristic zero to fields of
characteristic $p$, when the residual characteristic is large
enough, and vice versa. It should be emphasized that our statement
holds for quite general functions defined by integrals depending on
valued field variables. One should keep in mind that there is no
meaning in comparing values of individual  parameters in the
integrals or the integrals themselves between characteristic zero
and characteristic $p$. Our  transfer principle, which can be
considered as a wide generalization of the classical Ax-Kochen-Er{\v
s}ov result, should have a wide range of applications to $p$-adic
representation Theory and the Langlands Program, in particular to
various forms of the Fundamental Lemma. For instance, it applies
typically to the integrals occuring in  Jacquet-Ye conjecture
\cite{JY}, now a theorem thanks to work of Ng\^o \cite{Ngo} over
functions fields and Jacquet \cite{J} in general.

\medskip

Let us now review the content of the paper in more detail.
In section \ref{cef} we enlarge our Grothendieck rings in order to add exponentials.
In fact it is useful to consider not only exponentials of functions
with values in the valued field, but also exponentials of functions
with values in the residue field.
This is performed in a formal way by replacing the category
$\RDef_S$ considered in \cite{cl} - consisting of  certain objects
$X \rightarrow S$  - by a larger category
$\RDef_S^{\rm exp}$ consisting of the same
$X \rightarrow S$ together with functions
$g$ and $\xi$ on $X$ with values in the valued field, resp. the residue field.
We define a Grothendieck ring
$K_0 (\RDef_S^{\rm exp})$ generated by classes of objects
$(X, g, \xi)$ modulo  certain relations.
Here we have to add some new relations to the
classical ones already considered in \cite{cl}.
When $X \rightarrow S$ is the identity,
the class of $(X, g, 0)$, resp. $(X, 0, \xi)$,  corresponds to
the exponential of $g$, resp. the exponential of $\xi$.
One defines the ring $\cC (S)^{\rm exp}$ of motivic exponential functions
on $S$ by tensoring $K_0 (\RDef_S^{\rm exp})$ with the ring
$\cP (S)$ of
constructible
Presburger functions on  $S$. We are then able to state our main
results on integration of exponential functions in section
\ref{sec7e}. In particular we show that integrals with parameters of
functions in $\cC^{\rm exp}$ still lie in $\cC^{\rm exp}$. We first
directly construct integrals of exponential functions in relative
dimension 1 in section \ref{dim1} and then perform the general
construction in section \ref{sec8e}. As was the case in \cite{cl},
extensive use is made of the Denef-Pas cell decomposition Theorem.
Though some parts of our constructions and proofs are quite similar
to what we performed in \cite{cl}, or sometimes even folllow
directly from \cite{cl}, other require new ideas and additional work
specific to the exponential setting. As a first application, we
develop in section \ref{ft} the fundamentals  of
a
motivic Fourier transform. More precisely, there are two Fourier
transforms, the first one
over
residue field variables and the second one, which is more
interesting,  over valued field variables. Calculus with our valued
field Fourier transform is completely similar to the usual one.
Using convolution, we  define motivic Schwartz-Bruhat functions, and
we show that the valued field Fourier transform is involutive on
motivic Schwartz-Bruhat functions. We finally deduce Fourier
inversion for integrable functions with integrable Fourier
transform. In the following section \ref{sec:padics} we move to the
$p$-adic setting, defining the $p$-adic analogue of
 $\cC (S)^{\rm exp}$ and proving
stability under integration with parameters of these $p$-adic
constructible exponential functions. Such a result is the natural
extension to the exponential context of Denef's fundamental result
on stabillity of $p$-adic constructible functions under integration
with respect to parameters. This result of Denef greatly influenced
our work \cite{cl} and the present one. It has been later
generalized to the subanalytic case by the first author in
\cite{Ccell} and  \cite{Cexp}. In section \ref{transfert},  we loop
the loop by showing that motivic integration of constructible
exponential functions commutes with specialization to the
corresponding non archimedean ones, when the residue characteristic
is large enough. Finally, we end the paper by proving our
fundamental transfer Theorem, a form of which was already stated in
\cite{miami} when there is no exponential.
Let us note that in  their recent paper \cite{hk} Hrushovski and Kazhdan have also considered integrals of exponentials.

\medskip

Some of the results in this paper have been announced in
\cite{expnote}.

\section{Preliminaries}\label{prel}

\subsection{Definable subassignments and constructible functions}\label{dscf}
Let us start by recalling briefly some definitions and constructions from \cite{cl},
cf. also \cite{cr1},  \cite{cr2}.
We fix a field $k$ of characteristic
zero and we consider for any
field $K$ containing  $k$ the field of Laurent series
 $K \llp t \rrp$
endowed
with its natural valuation
\begin{equation}\ord : K \llp t \rrp^{\times} \longrightarrow \ZZ
\end{equation}
and with the angular component mapping
\begin{equation}\ac : K \llp t \rrp \longrightarrow K,
\end{equation}
defined by $\ac (x) = x t^{-\ord (x)} \mod t$ if $x \not=0$ and $\ac
(0) = 0$. We use the Denef-Pas language $\LPre$ which is a 3-sorted
language
\begin{equation}
(\LL_{\rm Val},\LL_{\rm Res},\LL_{\rm Ord},\ord,\ac)
\end{equation} with sorts
corresponding respectively to valued field, residue field and value group
variables.
The languages
$\LL_{\rm Val}$ and $\LL_{\rm Res}$ are equal to the ring language  $\LL_{\rm
Rings}=\{+,-,\cdot,0,1\}$, and for $\LL_{\rm Ord}$
we take the Presburger language
\begin{equation}
\LL_{\rm PR} = \{+, -, 0, 1, \leq\} \cup \{\equiv_n\ \mid n\in \NN,\
n
>1\},
\end{equation}
with  $\equiv_n$ the equivalence relation  modulo $n$.
Symbols $\ord$ and $\ac$
will be interpreted respectively as valuation and angular component,
so that
$(K \llp t \rrp,K,\ZZ)$ is a structure for $\LPre$.
We shall also add constants symbols in the
${\rm Val}$, resp. ${\rm Res}$, sort, for every element of $k \llp t \rrp$,
resp. $k$.

Let $\varphi$ be a  formula in the language $\LPre$
with respectively $m$, $n$ and $r$ free variables in the various sorts.
For every $K$ in $F_k$, the category of fields containing $k$,
we denote by $h_{\varphi} (K)$ the subset of
\begin{equation}
h [m, n, r] (K) := K \llp t \rrp^m \times K^n \times \ZZ^r
\end{equation} consisting of points
satisfying $\varphi$.
We call the assignment $K \mapsto h_{\varphi} (K)$ a definable
subassignment and we define a category $\Def_k$ whose objects
are
definable
subassignments.

More generally for $S$ in $\Def_k$, we denote by
$\Def_S$ the category of objects of $\Def_k$ over $S$.
We denote by
$\RDef_S$ the subcategory of
$\Def_S$ consisting of definable subassignments of
$S \times h [0, n, 0]$, for variable $n$, and
by $K_0 (\RDef_S)$ the corresponding Grothendieck ring.

We consider the ring
\begin{equation}A = \ZZ \Bigl[\LL, \LL^{-1}, \Bigl(\frac{1}{1 - \LL^{-i}}\Bigr)_{i > 0}\Bigr]
\end{equation}
and the subring $\cP (S)$ of the ring of functions from the set of
points of $S$ to $A$ generated by constant functions, definable
functions $S \rightarrow \ZZ$ and functions of the form
$\LL^{\beta}$ with $\beta$ definable $S \rightarrow \ZZ$. If $Y$
is a definable subassignment of  $S$, we denote by ${\bf 1}_Y$ the
function in  $\cP (S)$ with value  $1$ on  $Y$ and $0$ outside. We
denote by $\cP^0 (S)$ the subring of  $\cP (S) $ generated by such
functions and by  the constant function $\LL$. There is a morphism
$\cP^0 (S) \rightarrow K_0 (\RDef_S)$ sending ${\bf 1}_Y$ to the
class of $Y$  and sending $\LL $ to the class of $h [0, 1, 0]$.
Finally we
define
the ring of constructible motivic functions on $S$ by
\begin{equation}\cC (S) := K_0 (\RDef_S) \otimes_{\cP^0 (S)} \cP (S).
\end{equation}

To any algebraic subvariety $Z$ of
$\AA^m_{k  \llp t \rrp}$ we assign the definable subassignment
$h_Z$ of $h [m, 0, 0]$ given by  $h_Z (K) = Z (K \llp t \rrp)$.
The Zariski closure of a subassignment
$S$ of
$h [m, 0, 0]$ is the intersection  $W$ of all algebraic subvarieties
$Z$ of
$\AA^m_{k  \llp t \rrp}$ such that  $S \subset h_Z$. We set
$\dim S := \dim W$. More generally, if $S$ is a subassignment
of $h [m, n, r]$,
we define $\dim S$ to be
$\dim p (S)$ with $p$ the  projection
$h [m, n, r] \rightarrow h [m, 0, 0]$.
One proves, using results of \cite{Pas}
and \cite{vdDries}, that two isomorphic objects in $\Def_k$
have the same dimension.
For every non negative integer $d$,
we denote by  $\cC^{\leq d} (S)$ the ideal of  $\cC (S)$ generated by
the characteristic functions $\11_{Z}$ of definable subassignments
$Z$ of $S$ with $\dim Z \leq d$.
We set $C (S) = \oplus_d  C^d (S)$ with $C^d (S) := \cC^{\leq d} (S)
/ \cC^{\leq d-1} (S)$.

In \cite{cl}, we defined,
for $k$ a field of characteristic zero,
$S$ in $\Def_k$, and $Z$ in $\Def_S$,
a graded subgroup ${\rm I}_S C (Z) $ of
$C (Z) $ together with pushforward morphisms
\begin{equation}f_! : {\rm I}_S C (Z)
\longrightarrow
{\rm I}_S C (Y)
\end{equation} for every morphism $f : Z \rightarrow Y$ in $\Def_S$.
When $S$ is the final object $h [0, 0, 0]$ and $f$ is the morphism
$ Z \rightarrow S$, the morphism $f_!$ corresponds to motivic integration and we denote it by $\mu$.

Finally, fix $\Lambda$ in $\Def_k$.
Replacing dimension by relative dimension, we  defined
relative analogues $C (Z \rightarrow \Lambda)$ of
$C (Z)$ for $Z \rightarrow \Lambda$ in $\Def_{\Lambda}$ and
extended the above constructions to this relative
setting. In particular we constructed a morphism
\begin{equation}\mu_{\Lambda} :
{\rm I}_{\Lambda} C (Z \rightarrow \Lambda) \longrightarrow \cC (\Lambda)
= {\rm I}_{\Lambda} C (\Lambda \rightarrow \Lambda)\end{equation}
which corresponds  to motivic integration along  the fibers of the
morphism $Z \rightarrow \Lambda$.

\subsection{Cell decomposition}\label{ce}We now recall the definition of cells given in \cite{cl}, which is
a slight generalization of the one in  \cite{Pas}.

Let $C$ be a definable subassigment of $S$ where $S$ is in
$\Def_k$. Let $\alpha$, $\xi$, and $c$ be definable morphisms
$\alpha : C \rightarrow \ZZ$, $\xi : C \rightarrow h_{\GG_{m,
k}}$, and $c : C \rightarrow h[1,0,0]$. The cell $Z_{C, \alpha,
\xi, c}$ with basis $C$, order $\alpha$, center $c$, and angular
component $\xi$ is the definable subassignment of $S[1,0,0]$
defined by $y$ in $C$, $\ord (z - c (y)) = \alpha (y)$, and $\ac (z
- c(y)) = \xi (y)$, where $y$ lies in $S$  and $z$ in $h[1,0,0]$.
Similarly, if $c$ is a definable morphism $c : C \rightarrow
h[1,0,0]$, we define the cell $Z_{C, c}$ with center $c$ and basis
$C$ as the definable subassignment of $S[1,0,0]$ defined by $y \in
C$ and $z = c(y)$.

More generally, a definable subassignment $Z$ of $S[1,0,0]$ will
be called a 1-cell, resp.~a 0-cell, if there exists a definable
isomorphism
\begin{equation}\lambda : Z
\rightarrow Z_C = Z_{C, \alpha, \xi, c} \subset S[1,s,r],
\end{equation}
resp.~a definable isomorphism
\begin{equation}\lambda : Z
\rightarrow Z_C = Z_{C, c} \subset S [1,s,0],
\end{equation}
for some $s,r\geq 0$ and some $1$-cell $Z_{C, \alpha, \xi, c}$,
resp.~$0$-cell $Z_{C, c}$, such that the morphism
$\pi\circ\lambda$, with $\pi$ the projection on the
$S[1,0,0]$-factor, is the identity on $Z$.

We shall call the data $(\lambda, Z_{C, \alpha, \xi, c})$,
resp.~$(\lambda, Z_{C, c})$, sometimes written for short
$(\lambda, Z_C)$, a presentation of the cell $Z$.

One should note that $\lambda^{*}$ induces a canonical bijection
between $\cC (Z_C)$ and $\cC (Z)$.

In \cite{cl},
we proved
the  following
variant of the Denef-Pas Cell Decomposition Theorem \cite{Pas}:
\begin{theorem}\label{np}
Let $X$ be a definable subassignment  of $S[1,0,0]$ with $S$ in
$\Def_k$.
\begin{enumerate}
 \item[(1)]The subassigment $X$ is a finite disjoint union of cells.
 \item[(2)] For every $\varphi$ in $\cC (X)$ there exists a finite partition of $X$ into
cells $Z_i$ with presentation $(\lambda_i, Z_{C_i})$, such that
$\varphi_{| Z_i} = \lambda_i^{*} p_i^{*} (\psi_i)$, with $\psi_i$
in $\cC (C_i)$ and $p_i : Z_{C_i} \rightarrow C_i$ the projection.
Similar statements hold for $\varphi$ in $\cC (X)$, in $\cP(X)$,
and in $K_0(\RDef_X)$.
\end{enumerate}
\end{theorem}

We shall call a finite partition of $X$ into cells $Z_i$ as in
Theorem \ref{np} (1), resp.~\ref{np} (2) for a function $\varphi$, a
cell decomposition of $X$, resp.~a cell decomposition of $X$
adapted to $\varphi$.

The following result is already in  \cite{cl} (Theorem 7.5.3), except for (6) that
is new.

\begin{theorem}\label{isometryballs}Let $X$ be in $\Def_k$, $Z$ be a definable subassignment
of $X[1,0,0]$, and let $f:Z\to h[1,0,0]$ be a definable morphism.
There exists a cell decomposition of $Z$ into cells $Z_i$
such that the following conditions hold for every  $\xi$ in $C_i$,
for every  $K$ in $ {\rm Field}_{k(\xi)}$, and for every $1$-cell
$Z_i$ with presentation $\lambda_i:Z_i\to Z_{C_i}=Z_{C_i,\alpha_i,
\xi_i, c_i}$ and with projections $p_i:Z_{C_i}\to C_i$,
$\pi_i:Z_{C_i}\to h[1,0,0]$:
\begin{enumerate}
\item[(1)] The set $\pi_i(p_i^{-1}(\xi))(K)$ is either empty or a
ball of volume $\LL^{- \alpha_i(\xi)-1}$.  
 \item[(2)] When $\pi_i(p_i^{-1}(\xi))(K)$ is nonempty, the function
\begin{equation*}
g_{\xi,K} :
\begin{cases}
\pi_i(p_i^{-1}(\xi))(K)\to
K\llp t \rrp & \\
x\mapsto f\circ\lambda^{-1}(\xi,x) &
\end{cases}
\end{equation*}
is strictly analytic. 

\end{enumerate}
For each $i$ we can furthermore ensure that either $g_{\xi,K}$ is
constant or
\textup{(3)} up to \textup{(6)} hold.
\begin{enumerate}
 \item[(3)] There exists a definable morphism $\beta_i:C_i\to
h[0,0,1]$ such that
$$
\ord\, \frac{\partial }{\partial x}g_{\xi,K}(x)=\beta_i(\xi)
$$
for every $x$ in $\pi_i(p_i^{-1}(\xi))(K)$.  
 \item[(4)] When
$\pi_i(p_i^{-1}(\xi))(K)$ is nonempty, the map $g_{\xi,K}$ is a
bijection onto a ball of volume $\LL^{-\alpha_i(\xi) - 1 -
\beta_i(\xi)}$.

\item[(5)] For every
$x,y$ in $\pi_i(p_i^{-1}(\xi))(K)$, $\ord
(g_{\xi,K}(x)-g_{\xi,K}(y))=\beta_i(\xi)+\ord(x-y)$.

\item[(6)] There exists a morphism  $r_i:C_i\to h[1,0,0]$ such that
for every $x$ in $\pi_i(p_i^{-1}(\xi))(K)$
$$
g_{\xi,K}(x)=r_i(\xi)\ \mbox{ or }\
\ord(g_{\xi,K}(x)-r_i(\xi))\geq
\alpha_i(\xi) + \beta_i(\xi).
$$
\end{enumerate}
\end{theorem}

\begin{proof}
We only have to prove (6). First take a cell decomposition with
properties (1) up to (5). By replacing $X$ we may suppose that the
identity maps are presentations of the occurring cells. Then take a
cell decomposition of $W\subset X[1,0,0]$, with $W=p({\rm graph}f)$
and $p:X[2,0,0]\to X[1,0,0]$ the projection on the $X$ and $f$
coordinates. The centers of the so obtained cells are approximations
of $f$ as required in (6), and again by replacing $X$ one can assume
that the identity maps are presentations of the occurring cells. Now
take again a cell decomposition of $X$ such that properties (1) up
to (5) are fulfilled. Then automatically (6) is fulfilled as well.
\end{proof}

\section{Constructible exponential functions}\label{cef}

\subsection{Adding exponentials to Grothendieck rings}\label{grthe}

Let $Z$ be in $\Def_k$. We consider the category $\RDefe_{Z}$
whose objects are triples $(Y\to Z, \xi, g)$ with $Y$ in
$\RDef_{Z}$ and $\xi : Y \rightarrow h[0,1,0]$ and $g : Y
\rightarrow h[1,0,0]$ morphisms in $\Def_k$. A morphism $(Y'\to Z,
\xi', g') \rightarrow (Y\to Z, \xi, g)$ in $\RDefe_{Z}$ is a
morphism $h : Y' \rightarrow Y$ in $\Def_Z$ such that $\xi' = \xi
\circ h$ and $g' = g \circ h$. The functor sending $Y$ in
$\RDef_{Z}$ to $(Y, 0, 0)$, with $0$ denoting the constant
morphism with value $0$ in $h[0,1,0]$, resp.~$h[1, 0,0]$ being
fully faithful, we may consider $\RDef_{Z}$ as a full subcategory
of $\RDefe_{Z}$. We shall also consider the intermediate full
subcategory $\RDefw_{Z}$ consisting of objects $(Y, \xi, 0)$ with
$\xi : Y \rightarrow h[0,1,0]$ a morphism in $\Def_k$.

To the category $\RDefe_{Z}$ one assigns a Grothendieck ring $K_0
(\RDefe_{Z})$ defined as follows. As an abelian group it is the
quotient of the free abelian group over symbols $[Y \rightarrow Z,
\xi,g]$ with $(Y \rightarrow Z, \xi,g)$ in $\RDefe_{Z}$ by the
following four relations
\begin{equation}\label{r1}
[Y \rightarrow Z, \xi,g] = [Y' \rightarrow Z, \xi',g']
\end{equation}
for $(Y \rightarrow Z, \xi,g)$ isomorphic to $(Y' \rightarrow Z,
\xi',g')$,
\begin{equation}\label{r2}
\begin{split}
[(Y \cup Y') \rightarrow Z, \xi, g] & + [(Y \cap Y') \rightarrow Z,
\xi_{|Y \cap Y'},g_{|Y \cap Y'}]\\ &= [Y \rightarrow Z, \xi_{|Y},
g_{|Y}] + [Y' \rightarrow Z, \xi_{|Y'}, g_{|Y'}]
\end{split}
\end{equation}
for $Y$ and $Y'$ definable subassignments of some $X$ in $\RDef_Z$
and $\xi$, $g$ defined on $Y \cup Y'$,
 \begin{equation}\label{r3}
 [Y\to Z,\xi,g+h]=[Y\to Z,\xi+\bar h,g]
 \end{equation}
for $h:Y\to h[1,0,0]$ a definable morphism with $\ord(h(y))\geq 0$
for all $y$ in $Y$ and $\bar h$ the reduction of $h$ modulo $(t)$,
and
 \begin{equation}\label{r5}
 [Y[0,1,0]\to Z,\xi+p,g]=0
 \end{equation}
when $p:Y[0,1,0]\to h[0,1,0]$ is the projection and when
$Y[0,1,0]\to Z$, $g$, and $\xi$ factorize through the projection
$Y[0,1,0]\to Y$.

\begin{lem}\label{ring}
We may endow $K_0 (\RDefe_{Z})$ with a ring structure by setting
$$
[Y\to Z, \xi,g] \cdot  [Y'\to Z, \xi',g'] = [Y \otimes_Z Y'\to Z,
\xi\circ p_Y + \xi'\circ p_{Y'}, g\circ p_Y + g'\circ p_{Y'}],
$$
where $Y \otimes_Z Y'$ is the fiber product of $Y$ and $Y'$, $p_Y$
the projection to $Y$, and $p_{Y'}$ the projection to $Y'$.
\end{lem}
\begin{proof}
Clearly fiber product induces a commutative ring structure on the
free group on symbols $[Y \rightarrow Z, \xi,g]$ with $(Y
\rightarrow Z, \xi,g)$ in $\RDefe_{Z}$. The subgroup generated by
the four relations (\ref{r1}) up to (\ref{r5}) is an ideal of this
ring, hence, the quotient by this subgroup is a ring.
\end{proof}

Similarly, using relations (\ref{r1}), (\ref{r2}), (\ref{r5}), and
the subcategory $\RDefw_{Z}$, one may define the subring $K_0
(\RDefw_{Z})$ of $K_0 (\RDefe_{Z})$.

\subsubsection{Notation and abbreviations}
We write $\ee^\xi E(g) [Y\to Z]$ 
for $[Y\to Z, \xi,g]$. We abbreviate $\ee^{0}E(g)[Y\to Z]$ by
$E(g)[Y\to Z]$, $\ee^\xi E(0)[Y\to Z]$ by $\ee^\xi[Y\to Z]$, and
$\ee^{0}E(0)[Y\to Z]$ by $[Y\to Z]$. Likewise we write $\ee^\xi
E(g)$ for $\ee^\xi E(g)[Z\to Z]$, $E(g)$ for $\ee^0 E(g)[Z\to Z]$
and $\ee^\xi$ for $\ee^\xi E(0)[Z\to Z]$. Note that the element
$[Z\to Z]$ is the multiplicative unit of $K_0 (\RDefe_{Z})$.

\begin{lem}\label{l1e}
There are natural injections of rings $K_0 (\RDef_Z ) \rightarrow
K_0 (\RDefw_Z ) \rightarrow K_0 (\RDefe_{Z})$ sending $[Y\to Z]$ to
$[Y\to Z]$ and  $\ee^\xi[Y\to Z]$ to
 $\ee^\xi[Y\to Z]$. 
\end{lem}
\begin{proof}Both statements being similar, we prove that $i:K_0 (\RDef_Z ) \rightarrow
K_0 (\RDefw_Z )$ is an injection. Suppose that there are $a_1,a_2$ in
$K_0 (\RDef_Z )$ having the same image in $K_0 (\RDefw_Z )$. The
equality $i(a_1)=i(a_2)$ in $K_0 (\RDefw_Z )$ induces an equality in
the free group on symbols $[Y \rightarrow Z, \xi,0]$ of the form
\begin{equation}
\sum_i [Y_i \rightarrow Z, \xi_i,0] = \sum_j [Y_j \rightarrow Z,
\xi_j,0],
\end{equation}
by adding up relations. For each $i,j$, let $Y_i'\subset Y_i$ and
$Y_j'\subset Y_j$ be the subassignments defined by  $\xi_i=0$,
resp.~$\xi_j=0$. Then
\begin{equation}
\sum_i [Y'_i \rightarrow Z, 0,0] = \sum_j [Y'_j \rightarrow Z, 0,0]
\end{equation}
holds in the free group. Hence, $a_1=a_2$ in $K_0 (\RDef_Z )$.
\end{proof}

\subsection{Pull-back}For $f : Z \rightarrow Z'$ in $\Def_{k}$ we have a natural
pull-back morphism $f^* : K_0 (\RDefe_{Z'}) \rightarrow K_0
(\RDefe_{Z})$, induced by fiber product.

If $f : Z \rightarrow Z'$ is a morphism in  $\RDef_{Z'}$,
composition with $f$ induces a morphism $f_! : K_0 (\RDefe_{Z})
\rightarrow K_0 (\RDefe_{Z'}) $.

We have similar morphisms when we replace $\RDefe$ by $\RDefw$.

\subsection{Constructible exponential functions}\label{cstre} For
$Z$ in $\Def_k$ we define the ring $\cC  (Z)^{\rm exp}$ of
constructible exponential functions by
\begin{equation}
\cC  (Z)^{\rm exp}:=\cC (Z)\otimes_{K_0 (\RDef_{Z})} K_0
(\RDefe_{Z}) .
 \end{equation}

Note that the element $E({\rm id})$ of $\cC  (h[1,0,0])^{\rm
exp}$, with ${\rm id}$ the identity map on $h[1,0,0]$, can be seen
as a (non complex valued) additive character. Likewise, the
element $\ee^{\rm id}$ of $\cC (h[0,1,0])^{\rm exp}$, with ${\rm
id}$ the identity map on $h[0,1,0]$, can be seen as an additive
character on the residue field.

For every $d\geq 0$ we define $\cC ^{\leq d} (Z)^{\rm exp} $ as
the ideal of $\cC  (Z)^{\rm exp}$ generated by the characteristic
functions $\11_{Z'}$ of subassignments $Z'\subset Z$ of
dimension
 $\leq d$.

We set
\begin{equation}C (Z)^{\rm exp} = \oplus_d  C^d  (Z)^{\rm exp}
\end{equation}
 with
\begin{equation}
C^d  (Z)^{\rm exp} := \cC ^{\leq d} (Z)^{\rm exp} / \cC ^{\leq
d-1} (Z)^{\rm exp}.
\end{equation}
It is a graded abelian group, and also a $\cC  (Z)^{\rm
exp}$-module. We call elements of $C (Z)^{\rm exp}$ constructible
exponential Functions.

For $S$ in $\Def_k$ and $Z$ in $\Def_S$ we define the group $ {\rm
I}_S C (Z)^{\rm exp}$ of $S$-integrable constructible exponential
Functions by
 \begin{equation}
 {\rm I}_S C  (Z)^{\rm exp}:={\rm I}_S C
 (Z)\otimes_{K_0 (\RDef_{Z})} K_0 (\RDefe_{Z}) .
 \end{equation}
It is a graded subgroup of $C  (Z)^{\rm exp}$.

\begin{lem}\label{inj}
For every $Z$ in $\Def_k$, the natural morphisms of rings, resp.~of
graded groups, $\cC (Z)\to \cC (Z)^{\rm exp}$, $\cC ^{\leq d} (Z)
\rightarrow \cC ^{\leq d} (Z)^{\rm exp}$, resp.~$C (Z) \rightarrow C
(Z)^{\rm exp}$, ${\rm I}_S C (Z)\rightarrow {\rm I}_S C (Z)^{\rm
exp}$ are injective.
\end{lem}
\begin{proof}
This follows from Lemma \ref{l1e} by taking tensor products, and by
noting that $\cC ^{\leq d} (Z)^{\rm exp}$ is isomorphic to $
\cC^{\leq d} (Z)\otimes_{K_0 (\RDef_{Z})} K_0 (\RDefe_{Z})$.
\end{proof}

\begin{prop}\label{speprode}Let $S$ be in $\Def_k$ and let $W$ be a definable subassignment of $h [0, 0, m]$. The canonical
morphism
\begin{equation*}
K_0 (\RDefe_{S})  \otimes_{\cP^0  (S)}  \cP^0  (S \times W)
\longrightarrow K_0 (\RDefe_{S\times W})
\end{equation*}
is an isomorphism.
\end{prop}
\begin{proof}
This follows from the Denef-Pas quantifier elimination as stated in
\cite{cl}.
\end{proof}

\subsection{Inverse image of constructible exponential functions}\label{invime}
Let $f : Z \rightarrow Z'$ be a morphism in $\Def_k$. Since $f^*$ as
defined on $K_0 (\RDefe_{Z'})  $ and on $\cC(Z')$ is compatible with
the morphisms $K_0(\RDef_{Z'})\to \cC(Z')$ and $K_0(\RDef_{Z'})\to
K_0 (\RDefe_{Z'}) $, one gets by tensor product a natural pull-back
morphism $f^* : \cC (Z')^{\rm exp} \rightarrow \cC (Z)^{\rm exp}$.

\subsection{Push-forward for inclusions}\label{ince}
Let $i : Z \hookrightarrow Z'$ be the inclusion between two
definable subassignments $Z\subset Z'$. Extension by zero induces a
morphism $i_! : K_0 (\RDefe_{Z})  \rightarrow K_0 (\RDefe_{Z'})$.
Since this is compatible on $K_0(\RDef_Z)$ with $i_! : \cC (Z)
\rightarrow \cC (Z')$, we get, by tensor product, a morphism $i_! :
\cC (Z)^{\rm exp} \rightarrow \cC(Z')^{\rm exp}$. Because $i$ sends
subassignments of $Z$ to subassignments of $Z'$ of the same
dimension, there are group morphisms $i_! : \cC^{\leq d} (Z)^{\rm
exp} \rightarrow \cC^{\leq d} (Z')^{\rm exp}$, and graded group
morphisms $i_! : C (Z)^{\rm exp} \rightarrow C (Z')^{\rm exp}$. If
$Z'$ is in $\Def_S$ then $f_!$ clearly restricts to a morphism $f_!
:{\rm I}_S C (Z)^{\rm exp} \rightarrow {\rm I}_S C (Z')^{\rm exp}$.

\subsection{Push-forward for $k$-projections}\label{kproje}
Let $Y$ be in $\Def_k$ and let $Z$ be a definable subassignment of
$Y[0,r,0]$, for some $r\geq 0$. Denote by $f : Z \rightarrow Y$  the
morphism induced by projection. It follows from statement (1) in Proposition
5.2.1 of \cite{cl}  that the map $f_!: K_0 (\RDefe_{Z}) \rightarrow
K_0 (\RDefe_{Y}) $ induces a ring morphism $f_! : \cC (Z)^{\rm exp}
\rightarrow \cC (Y)^{\rm exp}$, and because $f$ sends subassignments
of $Z$ to subassignments of $Y$ of the same dimension, there are
group morphisms $f_! : \cC^{\leq d} (Z)^{\rm exp} \rightarrow
\cC^{\leq d} (Y)^{\rm exp}$, and graded group morphisms $f_! : C
(Z)^{\rm exp} \rightarrow C (Y)^{\rm exp}$. If $Y$ is in $\Def_S$
then $f_!$ clearly restricts to a morphism $f_! :{\rm I}_S C
(Z)^{\rm exp} \rightarrow {\rm I}_S C (Y)^{\rm exp}$. Note also that
the projection formula trivially holds in this setting, that is, for
every $\alpha$ in $\cC  (Y)^{\rm exp}$ and $\beta$ in $ C
(Z)^{\rm exp}$,
 $f_! (f^*(\alpha) \beta) = \alpha f_! (\beta)$.

\subsection{Push-forward for $\ZZ$-projections}\label{Zproje}
If $f : Z[0,0,m] \rightarrow Z$ is the projection and $Z$ is in
$\Def_S$, $m\geq 0$, then, by Proposition \ref{speprode} and by
the fact that $f$ preserves the dimension of definable
subassignments of $Z[0,0,m]$, the map
$f_!:{\rm I}_S C (Z[0,0,m]) \rightarrow {\rm I}_S C (Z)$

induces a graded group morphism $f_! : {\rm I}_S C (Z[0,0,m])^{\rm
exp} \rightarrow {\rm I}_S C (Z)^{\rm exp}$.

Lemma-Definition \ref{fub:kpres} below is a basic kind of Fubini Theorem
between the push forwards of \ref{kproje} and \ref{Zproje}, and
Lemma \ref{changekpres} is a basic form of the change of variables
formula.
\begin{def-lem}\label{fub:kpres}
Let $\varphi$ be in ${\rm I}_S C(Z[0,m,r])^{\rm exp}$ for some
$m,r\geq 0$ and some $Z$ in $\Def_S$ and let $f:Z[0,m,r]\to Z$ the
projection. Let $\pi_1,\ldots,\pi_{m+r}$ be any sequence of linear
projections of the form $Z[0,i,j]\to Z[0,i-1,j]$ or $Z[0,i,j]\to
Z[0,i,j-1]$ whose composition goes from $Z[0,m,r]$ to $Z$. Then,
$\pi_{m + r!}\circ\ldots\circ\pi_{1!}(\varphi)$ is independent of
the sequence $\pi_1,\ldots,\pi_{m+r}$ and we define $f_!(\varphi)$
to be this element.
\end{def-lem}
\begin{proof}
This follows from the fact that
\begin{equation*}
K_0 (\RDefe_{Z[0,m,0]}) \otimes_{K_0 (\RDefe_{Z}) } \cC (Z)^{\rm
exp} \otimes_{\cC (Z)}  \cC  (Z[0,0,r]) \longrightarrow \cC
(Z[0,m,r])^{\rm exp}
\end{equation*}
is an isomorphism.
\end{proof}

Let $\lambda :  S [0,n,r] \rightarrow S [0,n',r']$ be a  morphism
in $\Def_S$. Let $\varphi$ be a function in $\cC (S [0, n, r])^{\rm exp}$.
Assume $\varphi = \11_Z \varphi$ with $Z$ a
definable subassignment of $S [0, n, r]$ on which $\lambda$ is
injective. Thus $\lambda$ restricts to an isomorphism $\lambda'$
between $Z$ and $Z':= \lambda (Z)$. We define $\lambda_+(\varphi)$
in $\cC(S [0,n',r'])^{\rm exp}$ as $[i'_! (\lambda'{}^{-1})^* i^*
](\varphi)$, where $i $ and $i'$ denote respectively the
inclusions of $Z$ and $Z'$ in $S [0,n,r]$ and $S [0,n',r']$.

\begin{lem}\label{changekpres}
Let $\lambda : S [0, n, r] \rightarrow S [0, n' , r']$ be a
morphism in $\Def_S$. Let $\varphi$ be a function in $\cC(S [0, n,
r])^{\rm exp}$ such that  $\varphi = \11_Z \varphi$ with $Z$ a
definable subassignment of $S [0, n, r]$ on which $\lambda$ is
injective. Then $\varphi$ is in ${\rm I}_S C(S[0,n,r])^{\rm exp}$
if and only if $\lambda_+(\varphi)$ is in ${\rm I}_S
C(S[0,n',r'])^{\rm exp}$ and if this is the case then
$$
p_!(\varphi) = p'_!( \lambda_+ (\varphi)),
$$
with $p:S[0,n,r]\to S$ and $p':S[0,n',r']\to S$ the projections and
$p_!$ and $p'_!$ as in Lemma-Definition \ref{fub:kpres}.
\end{lem}
\begin{proof}
Consider the definable isomorphism
\begin{equation}
\lambda\times {\rm id}:S[0,n,r]\to S[0,n+n',r+r']
\end{equation}
with inverse $g$.
Since this is an isomorphism, $\varphi$ is $S$-integrable if and
only if
$g^*(\varphi)$
 is $S$-integrable.  By construction,
\begin{equation}
\pi'_!\left(g^*(\varphi)\right)=\lambda_+ (\varphi),
\end{equation}
\begin{equation}
\pi_!\left(g^*(\varphi)\right)=\varphi,
\end{equation}
with $\pi: S[0,n+n',r+r']\to S[0,n,r]$ and $\pi': S[0,n+n',r+r']\to
S[0,n',r']$ the projections. Now the Lemma follows from
Lemma-Definition \ref{fub:kpres}.
\end{proof}

\subsection{Relative setting}
Let us  fix  $\Lambda$ in $\Def_k$ that will play the role of a parameter space.
For  $Z$ in $\Def_{\Lambda}$,
we consider, similarly as in \cite{cl},
the ideal $\cC^{\leq d} (Z \rightarrow \Lambda)^{\rm exp}$
of $\cC (Z)^{\rm exp}$ generated by functions $\11_{Z'}$ with $Z'$ definable subassigment of $Z$
such that all fibers of $Z' \rightarrow \Lambda$ have dimension $\leq d$.
We set
\begin{equation}C (Z \rightarrow \Lambda)^{\rm exp} := \oplus_d  C^d (Z \rightarrow \Lambda)^{\rm exp}\end{equation}  with
\begin{equation}C^d (Z \rightarrow \Lambda) := \cC^{\leq d} (Z \rightarrow \Lambda)  / \cC^{\leq d-1} (Z \rightarrow \Lambda).\end{equation}
This graded abelian semigroup may be naturally identified with
\begin{equation}C (Z \rightarrow \Lambda) \otimes_{K_0 (\RDef_Z)} K_0 (\RDef_Z^{\rm exp}) .\end{equation}
For  $Z \rightarrow S$ a morphism in $\Def_{\Lambda}$, we set
\begin{equation}
{\rm
I}_S C (Z \rightarrow \Lambda)^{\rm exp} := {\rm I}_S C (Z
\rightarrow \Lambda) \otimes_{K_0 (\RDef_Z)} K_0 (\RDef_Z^{\rm exp})
.\end{equation}
Lemma \ref{inj} and all results and constructions in
(\ref{ince}), (\ref{kproje}), (\ref{Zproje}), including
Lemma-Definition \ref{fub:kpres} and Lemma
\ref{changekpres},  extend immediately with the same proofs to the relative setting.

\section{Integration of constructible exponential functions}\label{sec7e}

\subsection{The main result}We can now state
the
result on extending our construction of motivic integrals from
constructible functions to constructible exponential funcions.

\begin{theorem}\label{mte}Let $S$ be in
$\Def_k$. There is a unique functor from the category $\Def_S$ to
the category of abelian groups which sends  $Z$ to  ${\rm I}_S C
(Z)^{\rm exp}$, assigns to every morphism $f : Z \rightarrow Y$ in
$\Def_S$ a morphism $f_!  : {\rm I}_S C  (Z)^{\rm exp} \rightarrow
{\rm I}_S C  (Y)^{\rm exp}$ and which satisfies the following five
axioms:
\begin{enumerate}

\item[]{\textup{(A1) Compatibility:}} For every morphism $f : Z
\rightarrow Y$ in $\Def_S$, the map $f_!: {\rm I}_S C (Z)^{\rm
exp} \rightarrow {\rm I}_S C  (Y)^{\rm exp}$ is compatible with
the inclusions of groups ${\rm I}_S C  (Z) \to {\rm I}_S C
(Z)^{\rm exp}$ and ${\rm I}_S C  (Y) \to {\rm I}_S C (Y)^{\rm
exp}$ and with the map $f_!:{\rm I}_S C  (Z)\to {\rm I}_S C  (Y)$
as constructed in \cite{cl}.

\medskip

\item[]{\textup{(A2) Disjoint union: }}Let $Z$ and $Y$ be definable
subassignments in $\Def_S$. Assume $Z$, resp.~$Y$, is the disjoint
union of two definable subassignments $Z_1$ and $Z_2$, resp.~$Y_1$
and $Y_2$ of  $Z$, resp.~$Y$. Then, for every morphism $f : Z
\rightarrow Y$ in $\Def_S$, with $f (Z_i) \subset Y_i$ for $i = 1,
2$,
under the isomorphisms ${\rm I}_S C  (Z)^{\rm exp} \simeq {\rm
I}_S C (Z_1)^{\rm exp} \oplus {\rm I}_S C  (Z_2)^{\rm exp}$ and
${\rm I}_S C  (Y)^{\rm exp} \simeq {\rm I}_S C (Y_1)^{\rm exp}
\oplus {\rm I}_S C  (Y_2)^{\rm exp}$, we have $f_! = f_{1 !}
\oplus f_{2!}$, with $f_i:Z_i\to Y_i$ the restrictions of $f$.

\medskip

\item[]{\textup{(A3) Projection formula: }}
For every morphism $f : Z
\rightarrow Y$ in $\Def_S$, and every $\alpha$ in $\cC  (Y)^{\rm
exp}$ and $\beta$ in ${\rm I}_S C  (Z)^{\rm exp}$,
if $f^*(\alpha) \beta$ is in ${\rm I}_S C  (Z)^{\rm exp}$,
then $f_! (f^*(\alpha) \beta) = \alpha
f_! (\beta)$.

\medskip

\medskip
\item[]{\textup{(A4) Projection on $k$-variables: }}  Assume that
$f$ is the projection $f : Z = Y [0, n, 0] \rightarrow Y$ for some
$Y$ in $\Def_S$. For every  $\varphi$  in ${\rm I}_S C  (Z)^{\rm
exp}$,  $f_!(\varphi)$ is as constructed in \ref{kproje}.

\medskip
\item[]{\textup{(A5) Relative balls of large volume: }}Let $Y$ be in
$\Def_S$ and consider definable morphisms $\alpha : Y \rightarrow
\ZZ$, $\xi : Y \rightarrow h_{\GG_{m, k}}$, with $\GG_{m, k}$ the
multiplicative group $\AA^1_k \setminus \{0\}$. Suppose that $[\11_Z]$ is in ${\rm I}_S C (Z)^{\rm exp}$ and that
$Z$ is the definable subassignment of $Y [1, 0, 0]$ defined by
$\ord\, z = \alpha (y)$ and $\ac\, z = \xi (y)$, and $f : Z
\rightarrow Y$ is the morphism induced by the projection $Y[1,0,0]
\rightarrow Y$. If moreover $\alpha(y)<0$ holds for every $y$ in
$Y$, then
$$f_! (E(z)[\11_Z]) = 0.$$

\end{enumerate}
Moreover, these groups morphisms $f_!$ coincide with the group
morphisms constructed in \ref{Zproje} and in \ref{ince} in the
corresponding cases.
\end{theorem}

When $S = h [0, 0, 0]$, we write
${\rm I} C (Y)^{\rm exp}$
for
${\rm I}_S C (Y)^{\rm exp}$
and $\mu$ for the morphism $f_! : {\rm I} C (Y)^{\rm exp} \to {\rm I} C (h [0,0,0])^{\rm exp} = K_0 (\RDef_{h [0, 0, 0]}^{\rm exp}) \otimes_{\ZZ [\LL]} A$ when
$f : Y \rightarrow h [0, 0, 0]$ is the projection to the final object.

\subsection{Change of variables formula}We have the following analogue of
Theorem 12.1.1 of \cite{cl}.

\begin{theorem}[Change of variables formula]\label{cvfexp}Let $f : X \rightarrow Y$ be a definable
isomorphism
between definable subassignments of dimension $d$. Let
$\varphi$ be in $\cC^{\leq d} (Y)^{\rm exp}$ with a non zero class in  $C^{d} (Y)^{\rm exp}$. Then
$[f^{\ast} (\varphi)]$ belongs to  ${\rm I}_Y C^d (X)^{\rm exp}$ and
$$f_! ([f^{\ast} (\varphi)])
=
\LL^{\ordjac f \circ f^{-1}}[ \varphi].$$
\end{theorem}

\begin{proof}Similar to the proof of Theorem 12.1.1 of \cite{cl}.
It is enough to consider the cases where $f$ is an injection or a projection.
When $f$ is an injection the statement is true by construction. For projections,
one reduces to the case of the projection of a $0$-cell  as in Proposition
11.4.3 of \cite{cl}, which follows also by construction.
\end{proof}

\subsection{Relative version}Fix  $\Lambda$ in $\Def_k$.
The proof of Theorem \ref{mte} which we shall give in sections
\ref{dim1}
and
\ref{sec8e}
readily extends to the
following relative version:

\begin{theorem}\label{mter}Let  $\Lambda$ belong to  $\Def_k$ and let  $S$ belong to  $\Def_{\Lambda}$.
There exists a unique functor from $\Def_S$ to the category of abelian groups
assigning to any morphism  $f : Z \rightarrow Y$ in
$\Def_S$ a morphism $$f_{! \Lambda} :
{\rm I}_S C (Z\rightarrow \Lambda)^{\rm exp}
\rightarrow {\rm I}_S C (Y\longrightarrow \Lambda)^{\rm exp}$$
satisfying the analogues of axioms \textup{(A1)-(A5)} when replacing $ C (\_)$ by $C (\_ \rightarrow \Lambda)$. 
\end{theorem}

Let $Z$ be in $\Def_{\Lambda}$. For every point $\lambda$ of
$\Lambda$, we denote by $Z_{\lambda}$ the fiber of $Z$ at
$\lambda$, as defined in \cite{cl} 2.6. We have a natural restriction
morphism $i_{\lambda}^* : C (Z \rightarrow \Lambda)^{\rm exp} \rightarrow
C (Z_{\lambda})^{\rm exp}$,
which respects the grading.
Let $f : Z \rightarrow Y$
be a morphism in $\Def_{\Lambda}$ and
let $\varphi$ be  in
$C (Z \rightarrow \Lambda)^{\rm exp}$.
We denote by
$f_{\lambda} : Z_{\lambda}
\rightarrow Y_{\lambda}$ the restriction of $f$ to the fiber
$Z_{\lambda}$.
It follows from
Proposition 14.2.1 of \cite{cl} that
if $\varphi$ is in
${\rm I}_Y C (Z \rightarrow \Lambda)^{\rm exp}$,
then
$i_\lambda^*(\varphi)$
 is in
${\rm I}_{Y_{\lambda}} C (Z_{\lambda})^{\rm exp}$.
Furthermore, it follows from the constructions that,  then
\begin{equation}
i_{\lambda}^* (f_{!\Lambda}(\varphi)) =
f_{\lambda !} (i_{\lambda}^* (\varphi))
\end{equation}
 for every point
$\lambda$ of $\Lambda$.

When $S = \Lambda$ and $f$ is the morphism $Z \rightarrow \Lambda$,
we write $\mu_{\Lambda}$ for the morphism
$f_{!\Lambda} : {\rm I}_{\Lambda} C (Z \rightarrow \Lambda)^{\rm
exp} \rightarrow \cC (\Lambda)^{\rm exp} = {\rm I}_{\Lambda} C
(\Lambda \rightarrow \Lambda)^{\rm exp}$.

\begin{remark}\label{wifub}It follows from
the functoriality statement in Theorem \ref{mte}, resp. \ref{mter},  that for
$f : X \rightarrow Y$ and $g : Y \rightarrow Z$
in
$\Def_S$,
$(g \circ f)_{!} =
g_{!} \circ f_{!}$, resp.
$(g \circ f)_{!\Lambda} =
g_{!\Lambda} \circ f_{!\Lambda}$.
We shall sometimes refer to that property as ``Fubini Theorem''.
\end{remark}

\subsection{Global version}Once we have
Theorem \ref{cvfexp}
at our disposal it is possible to develop the theory on global subassignments, defined
by replacing affine spaces by general algebraic varieties, along the lines of
section 15 of \cite{cl}. Since this is essentially straightforward  we shall not give
more details here.

\section{Exponential integrals in dimension $1$}\label{dim1}

We shall start by constructing directly exponential integrals in relative
valued field
dimension $1$.
\subsection{}\label{4.1} Let $S$ be a definable
subassignment and consider a definable subassignment
$X\subset S[1,0,0]$ and denote by $\pi:X\to S$ the projection.
Let $M_X$ be the free group on symbols $[Y \rightarrow X,
\xi,g,\varphi]$ with $((Y \rightarrow X, \xi,g),\varphi)$ in
$\RDefe_{X}\times I_SC(X)$.

We construct a map
\begin{equation}
\pi_!:M_X\to C(S)^{\rm exp}
\end{equation}
and show that it factorizes through the natural surjective group

morphism
$M_X\to I_SC(X)^{\rm exp}$, thus obtaining a map
\begin{equation}
\pi_!:I_SC(X)^{\rm exp}\to C(S)^{\rm exp},
\end{equation}
which is the integral in relative dimension $1$.

Consider $a=[f:Y \rightarrow X, \xi,g,\varphi]$ in $M_X$.

We shall use a suitable isomorphism of the form $\lambda: Y\to
Y'\subset Y[0,n,r]$ which is an isomorphism over $Y$ and which is
adapted to $a$ in a certain sense.
Then we shall define $\pi_!$ by going through the commutative diagram
\begin{equation*}\xymatrix{
Y  \ar[d]_{\lambda}^{\cong} \ar[r]^{f} & X \ar[r]^{\pi} & S
\\
Y' \ar@{=}[d] \ar[r]^{\pi'} & S':=S[0,n+n_Y,r] \ar[ur]^{p}\\
A_1\cup A_2\cup B, & }
\end{equation*}
where $Y\subset X[0,n_Y,0]$, $\pi'$ and $p$ are the projections, and
where we will write $Y'$ as the disjoint union of $A_1,A_2$, and
$B$, along which $\pi'_!$ will be easy to define.

Let $a'$ be $[Y'\rightarrow Y', \xi':=\xi\circ\lambda^{-1},
g':=g\circ\lambda^{-1}, \varphi':=\lambda^{-1}{}^*f^*(\varphi)]$
in $M_{Y'}$. We construct $\lambda$, define $\pi'_!(a')$, and then
define $\pi_!(a)$ as $p_!\pi'_!(a')$, where $p_!$ is as in
Lemma-Definition \ref{fub:kpres}.

Write $\varphi'$ as $\sum_{i=1}^2[\varphi'_i]$ with $\varphi'_i$ in
$\cC^{\leq i}(Y)$ but not in $\cC^{\leq i-1}(Y)\setminus\{0\}$.
By Theorems \ref{np} and \ref{isometryballs} we can take an
appropriate $\lambda$ such that $\xi'=\tilde \xi\circ \pi'$ for some
$\tilde \xi:S'\to h[0,1,0]$, such that
$\varphi'_i=\pi'{}^*(\varphi_{i})$ for some $\varphi_{i}\in \cC
(S')$,
 and such that properties (1) up to (6) of Theorem
\ref{isometryballs} are fulfilled for $g'$. There are now uniquely
determined parts $A$, $B\subset Y'$, such that $g'(x,\cdot):y\mapsto
g'(x,y)$ is constant on $B_{x}$ for each $x$, and nonconstant and
injective on $A_{x}$ for each $x$, where $A_{x}=\{y\in h [1,0,0]\mid
(x,y)\in A\}$ and $B_{ x}=\{y\in h [1,0,0]\mid (x,y)\in B\}$ are the
fibers, and $x$ runs over $S'$.

It is clear that $g'_{|B}=\tilde g \circ \pi'_{|B}$ for a unique
definable $\tilde g: \pi'(B)\to h[1,0,0]$.

On the part $A$ we have to work differently.
By construction,
$A$ is a $1$-cell having the identity morphism as presentation. By
the previous use of Theorem \ref{isometryballs}, $A$ is the disjoint
union of $A_{1}$ and $A_{2}$, with
\begin{equation}
A_{1}:= \{(x,y)\in A\mid g'(x,\cdot)\mbox{ maps $A_{x}$ onto a
ball of volume $\LL^j$ with } j
\leq
 0\}.
\end{equation}
\begin{equation}
A_{2}:= \{(x,y)\in A\mid g'(x,\cdot)\mbox{ maps $A_{x}$ onto a
ball of volume $\LL^j$ with }j
>
 0\}.
\end{equation}
Note that the $A_{i}$ are cells which are their own presentation.

By property (6) of Theorem \ref{isometryballs}, there are
definable morphisms $r:S'\to h[1,0,0]$ and $\eta: S'\to h[0,1,0]$
such that
\begin{equation}
g'(x,y) - r(x) \equiv \eta(x)\bmod (t)
\end{equation}
for  $(x,y)\in A_{2}$, that is, either $ r_{i}(x) - g'(x,y)$ has
order $>0$ and $\eta(x)=0$, or, it has order $0$ and angular
component equal to $\eta(x)$.

\begin{def-lem}\label{int-1-cell}Consider  $\lambda$, $A_1$, $A_2$, $B$, $r$, and $\eta$
as constructed above. 
Define $\pi'_!(a')$ in $C (S')^{\rm exp}$ as
\begin{eqnarray*}
\pi'_!(a') & := &  \ee^{\tilde \xi}E(\tilde g) \,
\pi'_!\left(\11_{B}\varphi'\right)
 +
  \ee^{\tilde \xi + \eta}
E(r) \, \pi'_! \left(\11_{A_{2}}\varphi'\right)
\end{eqnarray*}
where $\pi'_!$ in the right hand side is as in \cite{cl}.
Then,
$\pi'_!(a')$ lies in ${\rm I}_S C (S')^{\rm exp}$ and is
independent of the
choice of $r$ and $\eta$.
Furthermore,
$p_{!}(\pi'_!(a'))$,
where $p_{!}$ is as in Lemma-Definition \ref{fub:kpres},
is
independent of the choice of $\lambda$,
so
we can define $\pi_!(a)$ in $C(S)^{\rm exp}$ as
$$\pi_!(a):= p_{!}(\pi'_!(a')).$$
We extend $\pi_!$ to a group
morphism $\pi_!:
M_X\to C(S)^{\rm exp}$.
\end{def-lem}
\begin{proof}
That $\pi'_!(\varphi')$ lies in ${\rm I}_S C (S')^{\rm exp}$ follows
easily from the fact that $\pi'_!\left(\11_{B}\varphi' \right)$ and
$\pi'_! \left(\11_{A_{1}}\varphi'\right)$ are in ${\rm I}_S C (S')$,
which is true by the main theorem of \cite{cl} and the fact that
$\varphi$ is in $I_SC(X)$. The independence from the choice of $r$
and $\eta$ is clear by relation (\ref{r3}) for $C(S)^{\rm exp}$.

We prove the independence from the choice of $\lambda:Y\to Y'$.
Although this is similar to the proof of Lemma-Definition 9.1.2 in
\cite{cl}, using furthermore relation (\ref{r5}),  let us give details.
If another map $\hat \lambda:Y\to \hat Y$ with the same properties
and with partition $\hat A_1,\hat A_2,\hat B$ is given,
 there exists a third map $\breve \lambda:Y\to \breve Y$ with the same properties and
with partition $\breve A_1, \breve A_2, \breve B$, such that $\breve
\lambda^{-1}(\breve B)$ contains both $\hat \lambda^{-1}(\hat B)$
and $\lambda^{-1}(B)$,
for example, the map $\breve \lambda:=\lambda\otimes_Y\hat \lambda$
has this property.
Necessarily, $\breve \lambda^{-1}(\breve B)$ is equal to the union
of $\lambda^{-1}(B)$ with a $0$-cell and is also equal to the union
of $\hat \lambda^{-1}(\hat B)$ with a $0$-cell, since $g'$ is
injective on $A$. It follows from  this,
from $\varphi'_i=\pi'{}^*(\varphi_{i})$
and since $A$ is a $1$-cell, that one has
\begin{equation}
p_!\pi'_!\left(\11_{B}\varphi'\right)=\hat
p_!\hat\pi'_!\left(\11_{\hat B}\hat\varphi'\right)=\breve
p_!\breve\pi'_!\left(\11_{\breve B}\breve\varphi'\right),
\end{equation}
with  obvious notation (this also follows from Lemma-Definition
9.1.2 in \cite{cl}). By Lemma-Definition \ref{fub:kpres}, one finds
\begin{equation}
p_!(\ee^{\tilde \xi}E(\tilde
g)\pi'_!\left(\11_{B}\varphi'\right))=\hat p_!( \ee^{\tilde{\hat
\xi}}E(\tilde{\hat g})\hat\pi'_!\left(\11_{\hat
B}\hat\varphi'\right))=\breve p_!(\ee^{\tilde{\breve
\xi}}E(\tilde{\breve g})
 \breve\pi'_!\left(\11_{\breve B}\breve\varphi'\right)).
\end{equation}

Next we compare the integrals over $A$, $\hat A$ and $\breve A$,
working still similarly as in the proof of Lemma-Definition 9.1.2 in \cite{cl}.

Note that automatically we have the following inclusions
\begin{equation}
\breve\lambda^{-1}(\breve A)\subset \lambda^{-1}(A),\qquad
\breve\lambda^{-1}(\breve A)\subset \hat\lambda^{-1}(\hat A),
\end{equation}
\begin{equation}
\breve\lambda^{-1}(\breve A_1)\subset \lambda^{-1}(A_1),\qquad
\breve\lambda^{-1}(\breve A_1)\subset \hat\lambda^{-1}(\hat A_1),
\end{equation}
but maybe not so for $A_2$. By Lemma-Definition 9.1.2 in \cite{cl}, one
has
\begin{equation}
p_!\pi'_! \left(\11_{A_{}}\varphi'\right)=\hat p_!\hat \pi'_!
\left(\11_{\hat A_{}}\hat \varphi'\right)=\breve p_!\breve
\pi'_! \left(\11_{\breve A_{}}\breve \varphi'\right).
\end{equation}

The subassignment $\breve\lambda^{-1}(\breve A_2)$ is equal to
$\lambda^{-1}(A_2)$ with a $1$-cell $C$ adjoined and with a $0$-cell
removed. In our construction, since
$\varphi'_i=\pi'{}^*(\varphi_{i})$,
the pushforward is stable under removing a $0$-cell from $A$ (or
from $A_2$). Relation (\ref{r5}) insures that the integral over $C$
is equal to zero, hence, the lemma is proved.
\end{proof}

\begin{def-lem}\label{f!1}
The map $\pi_!$ constructed in Lemma-Definition \ref{int-1-cell}
factorizes through the natural group homomorphism $M_X\to
I_SC(X)^{\rm exp}$. We write $\pi_!$ for the induced group
homomorphism
$$\pi_!:I_SC(X)^{\rm exp}\to C(S)^{\rm exp}.$$
\end{def-lem}
\begin{proof}
We have to check that $\pi_!:M_X\to C(S)^{\rm exp}$ factorizes
through the tensor product $N_X\otimes_{K_0(\RDef_X)}I_SC(X)$, with
$N_X$ the free group on symbols $[b]$ with $b$ in $\RDefe_{Z}$, and
then we have to check that it factorizes through the relations
(\ref{r1}) up to (\ref{r5}). By construction, $\pi_!:M_X\to
C(S)^{\rm exp}$ is bilinear in the factors $N_X$ and $I_SC(X)$ over
the ring $K_0(\RDef_X)$, hence it factorizes through
$N_X\otimes_{K_0(\RDef_X)}I_SC(X)$.

That $\pi_!$ factors through relation (\ref{r1}) is clear since its
definition is independent of the choice of $\lambda$,
cf.~Lemma-Definition \ref{int-1-cell}. Relation (\ref{r2}) is clear
by construction. Relation (\ref{r3}) also follows since we can
choose $\lambda$ in such a way that
$\bar{h}$
factors through the projection $\pi'$,
and then one can compare the original construction and definition
with the ones with $g$ replaced by $g+h$. We prove relation
(\ref{r5}). Assume that $a$ is of the form $[Y[0,1,0]\to
X,\xi+p,g,\varphi]$ with $p:Y[0,1,0]\to h[0,1,0]$ the projection and
that $Y[0,1,0]\to X$, $g$, and $\xi$ factorize through the
projection $Y[0,1,0]\to Y$. It follows by construction from
relation (\ref{r5}) that
$\pi_!(a)$ is zero.
\end{proof}

One deduces from Theorem \ref{isometryballs} the following change of
variables statement in relative dimension $1$:

\begin{prop}[Change of variables]\label{cv1e}
Let $X$ and $Y$ be  definable subassignments of dimension $r$ of
$S[1,0,0]$ for some $S$ in $\Def_k$ and let $f:X\to Y$ be a
definable isomorphism over $S$. Suppose that $X$ and $Y$ are
equidimensional of relative dimension $1$ relative to the projection
to $S$. Then, $\varphi$ is in ${\rm I}_S C^{r} (Y)^{\rm exp}$ if and
only if $\LL^{-\ordjac f}f^*(\varphi)$ is in ${\rm I}_S C^{r}
(X)^{\rm exp}$ and if this is the case then
$$
 \pi_{Y!}(\varphi)=\pi_{X!}(\LL^{-\ordjac f}f^*(\varphi))
$$
holds in
$C (S)^{\rm exp}$
with $\pi_Y:Y\to S$ and $\pi_X:X\to S$
the projections, $\pi_{Y!},\pi_{X!}$ as in Lemma-Definition
\ref{f!1}, and $\ordjac$ as in \cite{cl}.
\end{prop}
\begin{proof}
By linearity we may assume that $\varphi$ is of the form
\begin{equation}
\varphi = \ee^\xi E(g)[Z\to Y] \varphi_0,
\end{equation}
with $\varphi_0$ in ${\rm I}_S C (Y)$, $Z\subset Y[0,n,0]$, and
$\xi:Z\to h[0,1,0]$ and $g : Z \rightarrow h[1,0,0]$ definable
morphisms.
 By pulling back along $Z\to Y$, we may assume that $Z=Y$. Choose
$\lambda$ as in the construction of $\pi_{Y!}(\varphi)$ in
Lemma-Definition \ref{f!1}. By changing $\lambda$ we may suppose
that Theorem \ref{isometryballs} is also applied to the function
$p_1\circ f$, with $p_1:X\to h[1,0,0]$ the projection. But then
$\lambda\circ f$ can be used to compute $\pi_{X!}(\LL^{-\ordjac f}
f^*(\varphi))$ as in Lemma-Definition \ref{f!1} and is seen to be equal to
$\pi_{Y!}(\varphi)$.
\end{proof}

\section{Proof of Theorem \ref{mte}}\label{sec8e}
\subsection{Notation}If $p : X \rightarrow Z$ is a morphism
in $\RDef_Z$ and $\varphi$ a Function in
$C^i (Z)$ which is the class of $\psi$ in
$\cC^{\leq i} (Z)$, the class of $p^* (\psi)$ in
$C^i (X)$
depends only of $\varphi$, so we denote it by $p^* (\varphi)$.
This construction extends  by linearity to a morphism
$p^* : C (Z) \rightarrow C (X)$.

\subsection{}Replacing $K_0 (\RDefe_{Z}) $ by
the subring $K_0 (\RDefw_{Z}) $, one defines subobjects $\cC
(Z)^{\rm e}$, $C  (Z)^{\rm e}$ and ${\rm I}_S C  (Z)^{\rm e}$ of
$\cC  (Z)^{\rm exp}$, $C  (Z)^{\rm exp}$ and ${\rm I}_S C  (Z)^{\rm
exp}$ defined in \ref{cstre}, cf. Lemma \ref{l1e}.

Let us first prove that
Theorem \ref{mte} restricted to this setting hold:

\begin{prop}\label{mtetriv}Let $S$ be in
$\Def_k$. There is a unique functor from the category $\Def_S$ to
the category of abelian groups which sends  $Z$ to  ${\rm I}_S C
(Z)^{\rm e}$, assigns to every morphism $f : Z \rightarrow Y$ in
$\Def_S$ a morphism $f_!  : {\rm I}_S C  (Z)^{\rm e} \rightarrow
{\rm I}_S C  (Y)^{\rm e}$ and which satisfies
axioms \textup{(A1)} to \textup{(A4)} of Theorem \ref{mte}.
Moreover, these groups morphisms $f_!$ coincide with the group
morphisms constructed in \ref{Zproje} and in \ref{ince} in the
corresponding cases.
\end{prop}

\begin{proof}Let $f : Z \rightarrow Y$ be a morphism in
$\Def_S$. Consider
$\varphi$ in ${\rm I}_S C (Z)^{\rm e}$ of the form
\begin{equation}
{\bf e}^h  [X \rightarrow Z] \varphi_0
\end{equation}
with $p : X \rightarrow Z$ in $\RDef_Z$, $h : X \rightarrow h [0, 1, 0]$
and $ \varphi_0$ in ${\rm I}_S C (Z)$. We have
\begin{equation}
\varphi = p_! ({\bf e}^h p^*\varphi_0).
\end{equation}
Hence, if we denote by $\delta_{f,h} : X \rightarrow Y [0, 1, 0]$
the morphism
\begin{equation}x \longmapsto ((f \circ p) (x), h (x)),
\end{equation}
the axioms force to set
\begin{equation}\label{tpoz}
f_! (\varphi) :=
\pi_{Y!} ({\bf e}^{\xi} \delta_{f,h !} (p^* \varphi_0)),
\end{equation}
with $\pi_Y$ the projection $Y [0, 1, 0] \rightarrow Y$, $\xi$ the
canonical coordinate on the fibers of $\pi$,
and $\pi_{Y!}$ uniquely determined by (A4).
 Since ${\rm I}_S C (Z)^{\rm e}$ is generated by functions
$\varphi$ as above, this proves the unicity part of the statement.
For existence, one uses (\ref{tpoz}) to define $p_!$ by additivity.
Note that this definition is clearly compatible with the relations
involved in the definition of $C (Z)^{\rm e}$. Note also that (A2)
and (A1)
are obvious and that (A4), that is, compatibility  with
\ref{kproje}, is easily checked.
The projection formula (A3) follows easily from the projection
formula in \cite{cl}.

Now let us prove
functoriality, namely, that $g_! \circ f_! = (g \circ f)_!$ for
morphisms
 $f : Z \rightarrow Y$ and $g :
Y \rightarrow W$ in $\Def_S$.
 As above consider $\varphi$ in ${\rm
I}_S C (Z)^{\rm e}$ of the form ${\bf e}^h  [X \rightarrow Z]
\varphi_0= p_! ({\bf e}^h p^*\varphi_0)$ with $p : X \rightarrow Z$
in $\RDef_Z$,
$h : X \rightarrow h [0, 1, 0]$
and $ \varphi_0$ in ${\rm I}_S C (Z)$. We have
\begin{equation}
(g_! \circ f_!) (\varphi)=
g_! (\pi_{Y!} (({\bf e}^{\xi} \delta_{f, h !} (p^* {\varphi_0})))
\end{equation}
and
\begin{equation}
\begin{split}
 (g \circ f)_! (\varphi) & =
 \pi_{W !}({\bf e}^{\xi} \delta_{g \circ f, h !} (p^* {\varphi_0}))\\
 & = \pi_{W !}({\bf e}^{\xi} ((g \times {\rm id})_! \circ \delta_{f, h !}) (p^* {\varphi_0})),
 \end{split}
\end{equation}
hence it is enough to check that for every
$\psi$ in ${\rm I}_S C (Y [0, 1, 0])$
\begin{equation}\label{677}
g_! (\pi_{Y!} ({\bf e}^{\xi} \psi))
=
 \pi_{W !}({\bf e}^{\xi} ((g \times {\rm id})_! \psi)).
\end{equation}
We may assume $\psi$ is of the form $[p : X \rightarrow Y [0, 1, 0]] \, \pi_Y^* (\psi_0)$
with $\psi_0$ in ${\rm I}_S C (Y)$
the class of a function in some $\cC^i (Y)$, and with the above abuse of
notation.
Since
$\pi_{Y !} ({\bf e}^{\xi} \psi) = {\bf e}^{
\xi\circ p}
 [X \rightarrow Y] \, \psi_0$, we have
\begin{equation}
g_! (\pi_{Y!} ({\bf e}^{\xi} \psi))
=
 \pi_{W !}({\bf e}^{\xi} \delta_{g, \xi\circ p !} (
 p^* {\psi})).
\end{equation}
We now deduce
(\ref{677}) since
\begin{equation}
\begin{split}
\delta_{g, \xi \circ p!}
 (p^*
 {\psi} )
   & = (g \times {\rm id})_!    (  p_! ([X]p^* \pi_Y^* (\psi_0))))\\
&= (g \times {\rm id})_! ([X \rightarrow Y [0, 1, 0]] \, \pi_Y^*
(\psi_0))\\
& = (g \times {\rm id})_! (\psi).
\end{split}
\end{equation}
\end{proof}

\begin{remark}Note that in relative dimension 1, the morphisms $f_!$ in Proposition \ref{mtetriv}
coincide with those constructed in section \ref{dim1}.
\end{remark}

\subsection{Uniqueness}\label{munique}The proof is similar to the one in
Proposition \ref{mtetriv}. Let $f : Z \rightarrow Y$ be a morphism in
$\Def_S$. Consider
$\varphi$ in ${\rm I}_S C (Z)^{\rm exp}$ of the form
\begin{equation}
E (g) \ee^{h}   [X \rightarrow Z] \varphi_0
\end{equation}
with $p : X \rightarrow Z$ in $\RDef_Z$, $g : X \rightarrow h [1, 0,
0]$
and $h: X \rightarrow h [0, 1, 0]$ definable morphisms,
 and $
\varphi_0$ in ${\rm I}_S C (Z)$. We have
\begin{equation}
\varphi = p_! (E (g)
\ee^{h}
  p^*\varphi_0).
\end{equation}
Hence, if we denote by $\delta_{f,g,h} : X \rightarrow Y [1, 1, 0]$
the morphism sending $x$ to $((f \circ p) (x), g (x), h(x))$, the
axioms force to set
\begin{equation}\label{tppp}
f_! (\varphi) := \pi_{Y!} ( \pi_{Y[0,1,0]!} (E (x)
\ee^{\xi}  \delta_{f,g,h !}
   (p^* \varphi_0))),
\end{equation}
with $\pi_{Y[0,1,0]}:Y[1,1,0]\to Y[0,1,0] $ and $\pi_Y:Y [0, 1, 0]
\rightarrow Y$ the projections and $(x,\xi)$
the canonical coordinates on the fibers of $\pi$.
The map $\pi_{Y!}$ is determined by (A4) and for the map
$\pi_{Y[0,1,0]!}$ one is forced to use the construction of section
\ref{dim1}.

\subsection{Preliminaries}

Let $Z$ be in $\Def_S$. In Lemma-Definition \ref{fub:kpres}
we defined push-forwards morphisms
\begin{equation}
\pi_! : {\rm I}_S C (Z [0, r, m])^{\rm exp}
\longrightarrow
 {\rm I}_S C (Z)^{\rm exp}
\end{equation}
and
in section \ref{dim1} we constructed a pushforward
\begin{equation}
\pi_! : {\rm I}_S C (Z [1, 0, 0])^{\rm exp}
\longrightarrow
 {\rm I}_S C (Z)^{\rm exp},
\end{equation}
with $\pi$ denoting the projection.

We may mix these two constructions as follows:

\begin{def-lem}\label{uht}Let $Y$ be in $\Def_S$.
Let $\varphi$ be a Function in
${\rm I}_S C (Y [1, n, r])^{\rm exp}$.
Consider the following commutative diagram of projections
\begin{equation*}\label{thh}\xymatrix{
&Y [1, n, r] \ar[dl]_{\pi_1} \ar[dr]^{\pi'_1} \ar [dd]^{\pi} &\\
Y [1, 0, 0] \ar [dr]_{\pi_2}&&
Y [0, n, r] \ar [dl]^{\pi'_2}\\
&Y&.}
\end{equation*}
We have
$$\pi_{2 !} \pi_{1!}(\varphi) = \pi'_{2!} \pi'_{1!} (\varphi)$$
and we define $\pi_! (\varphi)$ to be the commun value
of
$\pi_{2 !} \pi_{1!}(\varphi)$ and
$ \pi'_{2!} \pi'_{1!} (\varphi)$.
\end{def-lem}

\begin{proof}The proof of Proposition-Definition 11.2.2  in \cite{cl} carries over to the present setting. \end{proof}

\subsection{Fubini for projections $Y [2, 0, 0] \rightarrow Y$}

\begin{prop}\label{ad}Let $Y$ be in $\Def_S$. Let $p : X \rightarrow Y [1, 0, 0]$ be in
$\RDef_{Y [1, 0, 0]}$ and $g : X \rightarrow h [1, 0, 0]$ be a morphism in $\Def_k$.
Denote by $\pi_Y$ the projection $Y [1, 0, 0] \rightarrow Y$ and
set $\gamma_g := (\pi_Y \circ p, g) : X \rightarrow Y [1, 0, 0]$.
For every Function $\psi$ in ${\rm I}_S C (Y [1, 0, 0])^{\rm e}$,
we have
$$
\pi_{Y !} (E (g)
[X\to Y]
\psi) = \pi_{Y !} (E (z) \gamma_{g !} (p^* (\psi))).
$$
\end{prop}

\begin{proof}
By using a construction with a cell decomposition adapted to $[X\to
Y[1,0,0]$, by pulling back, and by Lemma-Definition \ref{uht}, we
may assume $X= Y[1,0,0]$.
Similarly, by a cell decomposition construction using Theorem
\ref{isometryballs} and by Lemma-Definition \ref{uht}, we can reduce
to the case where $g$ is either constant or injective. When $g$ is
constant the statement is clear and when $g$ is injective it is a
direct consequence from Proposition \ref{cv1e}.
\end{proof}

Let $Y$ be in $\Def_S$. For $i = 1, 2$, we denote by $\pi_i : Y [2,
0, 0]\rightarrow Y [1, 0, 0]$ the projection $(y, z_1, z_2) \mapsto
(y, z_i)$ and  by $\pi_Y$ the projection $Y [1, 0, 0] \rightarrow
Y$.

\begin{prop}\label{spefub}Let $\psi$ be in
${\rm I}_S C (Y [2, 0, 0])^{\rm e}$.
Then
$$
\pi_{Y !} (\pi_{1 !} (E (z_2) \psi)) =
\pi_{Y !} (\pi_{2 !} (E (z_2) \psi)).
$$
\end{prop}

\begin{proof}We shall use bicells as defined in section 7.4 of \cite{cl}. By Proposition 7.4.1 of \cite{cl},
any definable subassignment $Z$ of $Y [2, 0, 0]$ admits a bicell
decomposition and, furthermore, for any $\varphi$ in $\cC (Z)$
there is a bicell decomposition of $Z$ adapted to $\varphi$ in the
sense of loc. cit. More generally, for any $\varphi$ in $\cC
(Z)^{\rm e}$, there is a bicell decomposition of $Z$ adapted to
$\varphi$. Indeed, this follows from the proof of loc. cit. and
the fact that statement (2) in Theorem \ref{np} still holds when
replacing $\cC (X)$ by $\cC (X)^{\rm e}$. Hence, we may assume
$\psi$ is the characteristic function of a bicell $Z$ in $Y [2, 0,
0]$. By the argument given at the beginning of the proof of
Proposition 11.2.4 of \cite{cl}, we may assume thet the bicell $Z$
is presented by the identity morphism.

Let us  consider first the case when $Z$ is a $(1, 1)$-bicell.
We start by the following special case:

\begin{lem}\label{easybicell}Let $C$ be a definable subassignment of
$Y$ and consider definable morphisms
$c : C \rightarrow h [1, 0, 0]$,
$\alpha$, $\beta : C \rightarrow h [0, 0, 1]$, and
$\xi$, $\eta : C \rightarrow h_{{\mathbf G}_{m, k}}$.
Consider the subassignment $Z$ of
$Y [2, 0, 0]$
defined by
\begin{gather*}
y \in C\\
\ord (z_1 - z_2) = \alpha (y)\\
\ac (z_1 - z_2) = \xi (y)\\
\ord (z_2 - c (y)) = \beta (y)\\
\ac (z_2- c(y)) = \eta (y).
\end{gather*}
Then
$$
\pi_{Y !} (\pi_{1 !} (E (z_2) \11_Z)) =
\pi_{Y !} (\pi_{2 !} (E (z_2) \11_Z)).
$$
\end{lem}

\begin{proof}[Proof of the Lemma]As in the proof of Lemma 11.2.5 of \cite{cl}, we may assume, after
partitioning $C$, that one of the following condition is satisfied everywhere on $C$:
\begin{enumerate}
\item [(1)]
 $\beta > \alpha$
 \item[(2)]$\beta < \alpha$
 \item[(3)]$\beta = \alpha$ and $\xi + \eta  \not=0$
 \item[(4)]$\beta = \alpha$ and $\xi + \eta  =0$.
\end{enumerate}
If condition (1) or (3) holds, $Z$ can be rewritten as a product of two $1$-cells, cf. loc. cit., and the result is clear.
If (2) is satisfied, then $Z$ is also defined by
\begin{gather*}
y \in C\\
\ord (z_1 - z_2) = \alpha (y)\\
\ac (z_1 - z_2) = \xi (y)\\
\ord (z_2 - c (y)) = \beta (y)\\
\ac (z_2- c(y)) = \eta (y).
\end{gather*}
and one computes
\begin{equation}
\pi_{Y !} (\pi_{1 !} (E (z_2) \11_Z))
=
E (c) I_{\alpha, \xi} I_{\beta, \eta}
\end{equation}
and
\begin{equation}
\pi_{Y !} (\pi_{2 !} (E (z_2) \11_Z))
=
E (c) \LL^{-\alpha - 1} I_{\beta, \eta}
\end{equation}
with $I_{\alpha, \xi}$, resp. $I_{\beta, \eta}$,
the integral of $E (z)$ over the subassigment of $ h [1, 0, 0]$ defined by
$\ord (z) = \alpha$ and $\ac (z) = \xi$, resp.
$\ord (z) = \beta$ and $\ac (z) = \eta$.
To deduce the requested equality note that $I_{\beta, \eta} = 0$ when $\beta <0$,
and that  when $\beta \geq 0$, $\alpha >0$, hence $I_{\alpha, \xi} = \LL^{-\alpha - 1}$.
The case of condition (4) also follows from
an  easy direct computation.
\end{proof}
When $Z$ is a $(1,1)$-bicell one proceeds similarly as in the proof of Proposition 11.5.4 of \cite{cl}.
More precidely assume  $Z$ is of the form
\begin{gather*}
y \in C\\
\ord (z_1 - d(y,z_2)) = \alpha (y)\\
\ac (z_1 - d (y,z_2)) = \xi (y)\\
\ord (z_2 - c (y)) = \beta (y)\\
\ac (z_2 - c(y)) = \eta (y).
\end{gather*}
If $d$ depends only on $y$, $Z$ is a product of $1$-cells and the
statement is clear. Otherwise $d(y,z_2)$
can be supposed injective, as in \cite{cl},
as function of $z_2$ for fixed $y$ and one deduces, with exactly the
same proof as in loc. cit., the statement from Lemma
\ref{easybicell} using the change of variables in relative dimension
$1$ (Proposition \ref{cv1e}).

When $Z$  is a $(1, 0)$-bicell one proceeds exactly as in loc. cit., using
change of variables in relative dimension $1$. In the  remaining cases of a
$(0,1)$ or a $(0,0)$-bicell, $Z$ is a product of cells, and the statement is clear.
\end{proof}

We may now prove the following version of Fubini Theorem:

\begin{prop}\label{fub}Let $\varphi$ be in
${\rm I}_S C (Y [2, 0, 0])^{\rm exp}$.
Then
$$
\pi_{Y !} (\pi_{1 !} (\varphi)) =
\pi_{Y !} (\pi_{2 !} (\varphi)).
$$
\end{prop}

\begin{proof}
By using a construction as in the proof of Lemma-Definition
\ref{uht}, but now with a bicell decomposition, we may assume
$\varphi = E (g)\psi$ with $g : Y [2, 0, 0] \rightarrow h [1, 0, 0]$
 in $\RDef_S$ and $\psi$ in ${\rm
I}_S C (Y [2, 0, 0])^{\rm e}$. By Proposition \ref{ad} we have
\begin{equation}
\begin{split}
\pi_{Y !} (\pi_{1_!}(\varphi)) &= \pi_{Y !} (\pi_{1_!}(E (g)
\psi))\\
&= \pi_{Y !} (\pi_{1_!}(E (z_2) \gamma_{1!} (\psi))),
\end{split}
\end{equation}
with
$\gamma_1 := (\pi_1, g) : X \rightarrow Y [2, 0, 0]$,
hence, by Proposition \ref{spefub}, we have
\begin{equation}
\begin{split}
\pi_{Y !} (\pi_{1_!}(\varphi)) &= \pi_{Y !} (\pi_{2_!}(E (z_2) \gamma_{1!}
 (\psi)))\\
&= \pi_{Y !} (E (z) (\pi_{2_!} \circ \gamma_{1!}) ( \psi))).
\end{split}
\end{equation}
The result follows, since
by Proposition \ref{mtetriv},
\begin{equation}
 (\pi_{2_!} \circ \gamma_{1!})
 (\psi) =
 (\pi_{2} \circ \gamma_{1})_!(\psi)
\end{equation}
and $\pi_2 \circ \gamma_1 : X \rightarrow Y [1, 0, 0]$ is the
morphism $(\pi,
 g)$,
with $\pi : Y [2, 0, 0] \rightarrow Y$ the projection, which is
independent from the order of the two variables in $h [2, 0, 0]$,
again, by Proposition \ref{mtetriv}.
\end{proof}

\subsection{Projections}

Let us  first consider the projection $p : Y [m, 0, 0] \rightarrow
Y$. Let $\varphi$ be in ${\rm I}_S C (Y [n, 0, 0])^{\rm exp}$. If $m
> 1$, we set, by induction on $m$, $p_! (\varphi) = \pi_{Y !}
(\pi_{1 !}(\varphi))$, where $\pi_{1} : Y [m, 0, 0] \rightarrow Y [m
- 1, 0, 0]$ is the projection on the first $m - 1$ coordinates. By
Proposition \ref{fub}, this definition is invariant under
permutation of coordinates.

In the general case of a projection $$p : Y [m, n, r]
\longrightarrow Y,$$ for any $\varphi$ in ${\rm I}_S C (Y [m, n,
r])^{\rm exp}$, we set $p_! (\varphi) := (p_{2!}\circ p_{1 !})
(\varphi)$ for
\begin{equation*}\xymatrix@1{
Y [m, n, r]  \ar[r]^>>>>{p_1}& Y [m, 0, 0] \ar[r]^>>>>{p_2} &Y.
}
\end{equation*}
It follows from Lemma-Definition \ref{fub:kpres},
Lemma-Definition \ref{uht} and Proposition \ref{fub}
that,
for any decomposition of $p$ into projections
\begin{equation*}\xymatrix@1{
Y [m, n, r]  \ar[r]^>>>>{p_1}& Y [m', n', r'] \ar[r]^>>>>{p_2} &Y,
}
\end{equation*}
 with $m' \leq m$, $n' \leq n$ and $r' \leq r$,
we have $p_! (\varphi) = (p_{2!}\circ p_{1 !}) (\varphi)$  and that
the definition of $p_! (\varphi)$  is invariant under permutation of
coordinates.

Now if $Z$ is a definable subassignment of some $h [m, n, r]$
and $\varphi$ is in ${\rm I}_S C (Y \times Z)^{\rm exp}$ we denote by
$\tilde \varphi$
the Function in
${\rm I}_S C (Y [m, n, r])^{\rm exp}$ which is obtained from
$\varphi$ by extension by  zero outside $Z$. If
$\varphi = \sum_{1 \leq i \leq j} \varphi_i  [X_i] E (g_i)$,
with $\varphi_i$ in ${\rm I}_S C (Y \times Z)^{\rm e}$,
$X_i$ in $\RDef_{Y \times Z}$     and
$g_i : X_i \rightarrow h[1, 0, 0]$,
we have, with a slight abuse of notations,
$\tilde \varphi = \sum_{1 \leq i \leq j} j_! (\varphi_i) j_! ([X_i]) E (g_i)$,
with $j : Z \rightarrow h [m, n, r]$ the inclusion.
Let us write the projection $p : Y \times Z \rightarrow Y$
as $\pi \circ j$, with $\pi$ the projection $Y [m, n, r] \rightarrow Y$,
and we set
\begin{equation}
p_! (\varphi) := \pi_! (\tilde \varphi)
\end{equation}
for $\varphi$  in ${\rm I}_S C (Y \times Z)^{\rm exp}$.

The projection formula (A3) trivially holds for $p_!$ and also
the following form of Fubini's Theorem.

\begin{prop}\label{hyptfu}Consider a diagram of projections
\begin{equation*}\xymatrix@1{
Y  \times Z \times W  \ar[r]^>>>>{p_1}& Y  \times Z \ar[r]^>>>>{p_2} &Y,
}
\end{equation*}
with $Z$ and $W$ in $\Def_k$.
Then for any
$\varphi$  in ${\rm I}_S C (Y \times Z \times W)^{\rm exp}$,
\begin{equation*}
(p_2 \circ p_1)_! (\varphi) = (p_{2 !} \circ p_{1!} ) (\varphi). \hfill  \qed
\end{equation*}
\end{prop}

\subsection{Definable injections}Let $i : X \rightarrow Y$ be a
morphism in $\Def_S$ which is a definable injection and let $g : X
\rightarrow h [1, 0, 0]$ and
$h:X \rightarrow h [0, 1, 0]$
be morphisms. Now if
$\varphi = \varphi_0 E(g)\ee^h [W\to X] $
 lies in ${\rm I}_S C (X)^{\rm exp}$
with $ \varphi_0$ in ${\rm I}_S C (X)$, $W$ in $\RDef_{X}$,
and if we write $W'$ for the unique element of $\RDef_Y$ such that
for each $x\in X$ the fiber $W_x$ equals $W'_{i(x)}$, we set
\begin{equation}
i_! (\varphi) := i_! (\varphi_0) E(g_{W'})\ee^{h_{W'}} [W'\to Y],
\end{equation}
where $h_{W'}$ and $g_{W'}$ are $h$ and $g$ seen on $W'$.

This definition extends uniquely
by linearity to give a morphism $i_! : {\rm I}_S C (X)^{\rm exp}
\rightarrow {\rm I}_S C (Y)^{\rm exp}$. Also it is quite clear that
if $j : Y \rightarrow Z$ is another definable injection in $\Def_S$,
$(j \circ i)_! =j_! \circ i_!$.

\subsection{General case}To define $f_!$ for $f : X \rightarrow Y$ a general morphism
in $\Def_S$, we proceed as in \cite{cl}. We decompose $f$ as a composition
$f = \pi_f \circ i_f$ with $i_f$ the definable injection
$X \rightarrow X \times Y$ given by
$x \mapsto (x, f(x))$
and
$\pi_f : X \times Y \rightarrow Y$ the canonical projection, and we set
$f_! = \pi_{f!}  \circ i_{f !}$.

It is quite clear  that when $f$ is an injection the new definition coincides with the previous one.
Also, when $f$ is a projection $Y \times Z \rightarrow Y$, the new definition
coincides with the previous one.
Indeed, the analogue of Lemma 11.5.2 of \cite{cl} holds for similar reasons
in the present setting, so the proof in Lemma 11.6.1 of
\cite{cl} extends directly.

\subsection{End of proof}We still have to check that
if $f : X \rightarrow Y$ and
$g: Y \rightarrow Z$ are morphisms in $\Def_S$, then
$g_! \circ f_! = (g \circ f)_!$. This is proved in a formal way exactly as in
Proposition 11.7.1 of \cite{cl}, since the analogue of Lemma 11.7.2 in \cite{cl}
holds in the present setting. Axioms (A1)-(A5) follow directly by construction.

\section{Fourier transformation}\label{ft}

Let $p : X \rightarrow \Lambda$ be a morphism in $\Def_k$, with all
fibers of relative
dimension
$d$. We shall denote by $\cI_{\Lambda} (X)^{\rm exp}$ or $\cI_{p}
(X)^{\rm exp}$ the $\cC (\Lambda)^{\rm exp}$-module of functions
$\varphi$ in $\cC (X)^{\rm exp}$ whose class $[\varphi]$ in $C^d (X
\rightarrow \Lambda)^{\rm exp}$ lies in ${\rm I}_{\Lambda} C ( Z
\rightarrow \Lambda)^{\rm exp}$. We shall also write $\mu_{\Lambda}
(\varphi)$ or $\mu_{p} (\varphi)$ to denote
the function
$\mu_{\Lambda} ([\varphi])$
in
$\cC (\Lambda)^{\rm exp}$.

\subsection{Fourier transformation over the residue field}Fix $\Lambda$ in $\Def_k$ and  an integer $d \geq 0$.
We consider the subassignment $\Lambda [0, 2d, 0]$ with first $d$
residue field coordinates $x = (x_1, \dots, x_d)$ and last $d$
residue coordinates $y = (y_1, \dots, y_d)$ and denote by $p_1 :
\Lambda [0, 2d, 0] \rightarrow \Lambda [0, d, 0]$ and $p_2 : \Lambda
[0, 2d, 0] \rightarrow \Lambda [0, d, 0]$ the projection onto  the
$x$-variables and $y$-variables, respectively. We consider the
function $\mathbf{e} (xy) := \mathbf{e} (\sum_{1 \leq i \leq d}x_i
y_i)$ in $\cC (\Lambda [0, 2d, 0])^{\rm exp}$.

We define the Fourier transform
\begin{equation}
\Fr : \cC (\Lambda [0, d, 0])^{\rm exp}
\longrightarrow
\cC (\Lambda [0, d, 0])^{\rm exp}
\end{equation}
by
\begin{equation}
\Fr (\varphi) := p_{1 ! \Lambda [0, d, 0]} ([\mathbf{e} (xy) p_2^* (\varphi)]),
\end{equation}
for
$\varphi$
in $\cC (\Lambda [0, d, 0])^{\rm exp}$.
The morphism
$\Fr$ is $\cC (\Lambda)^{\rm exp}$-linear.

For $\varphi$ in
$\cC (\Lambda [0, d, 0])^{\rm exp}$ we write
$\check{\varphi}$ for $\iota^* (\varphi)$, with $\iota : \Lambda [0, d, 0]
\rightarrow \Lambda [0, d, 0]$ the $\Lambda$-morphism sending
$x$ to $- x$.

\begin{theorem}\label{frr}Let $\varphi$ be in $ \cC (\Lambda [0, d, 0])^{\rm exp}$. We have
$$
(\Fr \circ \Fr) (\varphi)
= \LL^d \,  \check{\varphi}.
$$
\end{theorem}

\begin{proof}
The proof is essentially the same as the standard one for finite fields.
More precisely,
we  work on $\Lambda [0, 3d, 0]$
with coordinates $(x_1, \dots, x_d, y_1, \dots, y_d, z_1, \dots, z_d)$.
We shall denote by $x : \Lambda [0, 3d, 0] \rightarrow \Lambda [0, d, 0]$,
$(x, y )  : \Lambda [0, 3d, 0] \rightarrow \Lambda [0, 2d, 0]$
the projections onto  the corresponding components, etc.
If $f$ is a function on $ \Lambda [0, d, 0]$ we shall write
$f (x)$ instead of $x^* f$, etc.
By induction and Fubini Theorem (cf. Remark \ref{wifub}) we may assume $d = 1$.

Let $\varphi$ be in $ \cC (\Lambda [0, d, 0])^{\rm exp}$. We have
\begin{equation}
\begin{split}
(\Fr \circ \Fr) (\varphi)
&= \mu_x [\mathbf{e} (y (x + z)) \varphi (z)] \\
&=\mu_x [\mathbf{e} (y u ) \varphi ( u - x)],
\end{split}
\end{equation}
after performing the change of variables
$(x, y, z ) \mapsto (x, y, u = x + z)$.
Since
\begin{equation}
\LL \,  \check{\varphi} (x) = \mu_x [\11_{u = 0}\11_{y = y}\varphi (
- x)],
\end{equation}
we have
\begin{equation}
((\Fr \circ \Fr) (\varphi) -
\LL \,  \check{\varphi}) (x) =
\mu_x [\11_{u \not=0}\mathbf{e} (y u ) \varphi ( u - x)],
\end{equation}
which, after performing the change of variables
$(x, y, u ) \mapsto (x, w = yu, u)$, may be rewritten as
\begin{equation}
((\Fr \circ \Fr) (\varphi) -
\LL \,  \check{\varphi}) (x) =
\mu_x [\11_{u \not=0}\mathbf{e} (w ) \varphi ( u - x)],
\end{equation}
whose right hand side  is zero by Fubini Theorem and relation (\ref{r5}).
\end{proof}

\subsection{Fourier transformation over the valued field}\label{sfv}Fix $\Lambda$ in $\Def_k$ and  an integer $d \geq 0$.
We consider the subassignment $\Lambda [2d, 0, 0]$ with first $d$
valued field coordinates $x = (x_1, \dots, x_d)$ and last $d$ valued
field coordinates $y = (y_1, \dots, y_d)$ and denote by $p_1 :
\Lambda [2d, 0, 0] \rightarrow \Lambda [d, 0, 0]$ and $p_2 : \Lambda
[2d, 0, 0] \rightarrow \Lambda [d, 0, 0]$ the projection onto  the
$x$-variables and $y$-variables, respectively. We consider the
function $E (xy) := E (\sum_{1 \leq i \leq d}x_i y_i)$ in $\cC
(\Lambda [2d, 0, 0])^{\rm exp}$.

\begin{lem}\label{trivi}Let $\varphi$ be
in $\cI_{\Lambda} (\Lambda [d, 0, 0])^{\rm exp}$. The class $[E
(xy) p_2^* (\varphi)]$ of $E (xy) p_2^* (\varphi)$ in $C^d (p_1 :
\Lambda [2d, 0, 0] \rightarrow \Lambda [d, 0, 0])^{\rm exp}$ is
integrable rel.~$p_1$.
\end{lem}

\begin{proof}Indeed, it follows by construction that,
if  $\psi$ is a function in $\cC (\Lambda [d, 0, 0])$ whose class $[\psi]$
in
$C^d (\Lambda [d, 0, 0] \rightarrow \Lambda)$ is
$\Lambda$-integrable, the class of the pull-back
$p_2^* (\psi)$
in
$C^d (p_1 :
\Lambda [2d, 0, 0] \rightarrow \Lambda [d, 0, 0])$ is
integrable rel.~$p_1$.
The statement follows.
\end{proof}

Thanks to Lemma \ref{trivi},
one may define the Fourier transformation
\begin{equation}
\Fv : \cI_{\Lambda} (\Lambda [d, 0, 0])^{\rm exp}
\longrightarrow
\cC (\Lambda [d, 0, 0])^{\rm exp}
\end{equation}
by
\begin{equation}
\Fv (\varphi) := \mu_{p_1} (E (xy) p_2^* (\varphi)),
\end{equation}
for
$\varphi$
in $\cI_{\Lambda} (\Lambda [d, 0, 0])^{\rm exp}$.
The morphism
$\Fv$ is $\cC (\Lambda)^{\rm exp}$-linear.

\subsection{}\label{ex}Let us compute some simple examples.
Consider definable functions $\alpha : \Lambda \rightarrow \ZZ$  and
$\xi = (\xi_1, \dots, \xi_d) : \Lambda \rightarrow h [0, d, 0]$ with
$\xi_i$ nowhere zero
and set
$$Z_{\alpha} := \{(\lambda, x = (x_1, \dots, x_d)) \in \Lambda [d, 0, 0] \mid \ord (x_i) \geq \alpha (\lambda)\},$$
$$W_{\alpha} := \{(\lambda, x= (x_1, \dots, x_d)) \in \Lambda [d, 0, 0] \mid \ord (x_i) = \alpha (\lambda)\},$$
$$W_{\alpha, \xi } := \{(\lambda, x= (x_1, \dots, x_d)) \in \Lambda [d, 0, 0] \mid \ord (x_i) = \alpha (\lambda),
\ac (x_i) = \xi_i (\lambda)\},$$
$$\varphi_{\alpha} :={\bf{1}}_{Z_{\alpha}}, \quad
\psi_{\alpha} :={\bf{1}}_{W_{\alpha}}, \quad \rm{and}\quad
\psi_{\alpha, \xi} :={\bf{1}}_{W_{\alpha, \xi}}.$$

\begin{prop}\label{form}The following formulas hold:
\begin{enumerate}
\item[(1)] $\Fv (\varphi_{\alpha}) = \LL^{- d \alpha}
\varphi_{-\alpha + 1}$.
\item[(2)] $\Fv (\psi_{\alpha}) =
\LL^{- d \alpha}  \varphi_{-\alpha + 1}
-  \LL^{- d\alpha - d}  \varphi_{-\alpha}$.
\item[(3)]
$\Fv (\psi_{\alpha, \xi}) = \LL^{- d\alpha - d} (\varphi_{- \alpha +
1}
 +
  \ee(i)\psi_{-\alpha}),$ with $i$ the morphism $\Lambda [d, 0, 0] \rightarrow \Lambda
[d, 1, 0]$ given by $(\lambda, x) \mapsto (\lambda, x, \sum_i
\xi_i(\lambda) \ac\, x_i)$.
\end{enumerate}
\end{prop}

\begin{proof}By induction on $d$, we may assume $d = 1$.
Let us start by proving (3). It is enough to check that the
restriction of $\Fv (\psi_{\alpha, \xi})$ to the subassigment
defined by $\ord\, x = \beta$ and $\ac\, x = \eta$ is equal to $0$
if $\alpha + \beta <0$, to $\LL^{- \alpha - 1}$ if $\alpha + \beta
>0$, and is equal to ${\bf e} (\xi \eta) \LL^{- \alpha - 1}$ if
$\alpha + \beta =0$.
The case $\alpha + \beta <0$ follows from (A5) of Theorem
\ref{mte}. The cases  $\alpha + \beta >0$  and $\alpha + \beta =0$
follow from relation (\ref{r3}) and the construction in \cite{cl}.
Cases (1) and (2) are more easy. The reader may also choose to
prove first the case of $\alpha =0$ and deduce the case of general
$\alpha$ from it.
\end{proof}

We calculate some more examples.

\begin{lem}\label{hypertriv} Assume $d = 1$ and let $\gamma :  \Lambda \rightarrow \ZZ$
and $\xi : \Lambda \rightarrow h [0, 1, 0]$
be definable functions. Then
\begin{enumerate}
\item[(1)] If $\gamma > 0$ on $\Lambda$,
$\mu_{\Lambda} (\psi_{\gamma, \xi } E (x)) = \LL^{- \gamma - 1}$,
$\mu_{\Lambda} (\psi_{\gamma} E (x)) = (\LL - 1) \LL^{- \gamma - 1}$
and
$\mu_{\Lambda} (\varphi_{\gamma} E (x)) =  \LL^{- \gamma }$.
\item[(2)]
If $\gamma < 0$ on $\Lambda$,
$\mu_{\Lambda} (\psi_{\gamma, \xi } E (x)) = \mu_{\Lambda} (\psi_{\gamma} E (x)) = 0$.
\item[(3)] If $\gamma = 0$ on $\Lambda$,
$\mu_{\Lambda} (\psi_{\gamma, \xi} E (x)) = {\bf e} (\xi) \LL^{- 1}$
and $\mu_{\Lambda} (\psi_{\gamma} E (x)) = - \LL^{- 1}$.
\end{enumerate}
\end{lem}

\begin{proof}Statement (1) and the first part of (3) are  obvious
from relation (\ref{r3}) and the construction in \cite{cl},
and (2) follows from (A5) of Theorem \ref{mte}.

The last part of (3) follows from the first part using cell
decomposition, since, by relation (\ref{r5}),  $\mu_{\Lambda} (i_!
i^* {\bf e} (\xi)) = - 1$, with $\xi$ the residue field variable
on $\Lambda [0, 1, 0]$ and $i$ the inclusion of the subassignment
defined by $\xi \not= 0$ in $\Lambda [0, 1, 0]$.
\end{proof}

If follows readily from Proposition \ref{form}
that
\begin{equation}\label{ef1}
\Fv \circ \Fv (\varphi_{\alpha}) =  \LL^{- d} \, \varphi_{\alpha}
\end{equation}
and
\begin{equation}
\Fv \circ \Fv (\psi_{\alpha}) =   \LL^{- d} \, \psi_{\alpha}.
\end{equation}
The corresponding statement for $\psi_{\alpha, \xi}$ will follow
from the general Fourier inversion for Schwartz-Bruhat functions
we shall prove  in Theorem \ref{sinv}. Though not at all
necessary, let us provide a direct proof of that fact:

\begin{prop}\label{gty}The following holds
$$\Fv \circ \Fv (\psi_{\alpha, \xi}) = \LL^{- d} \, \psi_{\alpha, - \xi}.$$
\end{prop}

\begin{proof}We may assume $d = 1$.
By (1) and (3) of Proposition \ref{form}, it is enough to
calculate $\Fv (\ee(i)\psi_{-\alpha}),$ with $i$ the morphism
$\Lambda [1, 0, 0] \rightarrow \Lambda [1, 1, 0]$ given by
$(\lambda, x) \mapsto (\lambda, x, \xi(\lambda) \ac\, x)$. But
this can be done similarly as (3) of Proposition \ref{form} and
the last part of (3) of Lemma \ref{hypertriv}.
\end{proof}

\subsection{Convolution}
We denote by $x + y$ the morphism $\Lambda [2d, 0, 0] \rightarrow
\Lambda [d, 0, 0]$ given by $(x_1, \dots, x_d, y_1, \dots, y_d)
\mapsto (x_1 + y_1, \dots, x_d + y_d)$. We shall also work on
$\Lambda [3d, 0, 0]$ with coordinates $(x_1, \dots, x_d, y_1,
\dots, y_d, z_1, \dots, z_d)$. We shall denote by $x : \Lambda
[3d, 0, 0] \rightarrow \Lambda [d, 0, 0]$, $(x, y )  : \Lambda
[3d, 0, 0] \rightarrow \Lambda [2d, 0, 0]$ the projections onto
the corresponding components, etc. If $f$ is a function on $
\Lambda [d, 0, 0]$ we shall write $f (x)$ instead of $x^* f$, etc.

\begin{def-prop}Let $f$ and $g$ be two functions in
$\cI_{\Lambda} (\Lambda [d, 0, 0])^{\rm exp}$.
The function
$p_1^* (f) p_2^* (g)$
lies in
$\cI_{x + y} (\Lambda [2d, 0, 0])^{\rm exp}$ and
the function
\begin{equation*}
f \ast g := \mu_{x + y}(p_1^* (f) p_2^* (g))
\end{equation*}
lies in $\cI_{\Lambda} (\Lambda [d, 0, 0])^{\rm exp}$,
where $p_1$ and $p_2$ are as in section \ref{sfv}.
\end{def-prop}

\begin{proof}It follows directly from Theorem 14.1.1 in \cite{cl}
that, if
$\varphi$ and $\psi$ are functions in $\cC (\Lambda [d, 0, 0])$
whose class in $C^d (\Lambda [d, 0, 0] \rightarrow \Lambda)$ are in
${\rm I}_{\Lambda} C (\Lambda [d, 0, 0] \rightarrow \Lambda)$,
then the class of $p_1^* (\varphi) p_2^* (\psi)$ in
$C^d (\Lambda [2d, 0, 0] \rightarrow \Lambda)$ lies in
${\rm I}_{\Lambda} C (\Lambda [2d, 0, 0] \rightarrow \Lambda)$.
One deduces that
the function
$p_1^* (f) p_2^* (g)$
lies in $\cI_{\Lambda} (\Lambda [2d, 0, 0])^{\rm exp}$, hence also in
$\cI_{x + y} (\Lambda [2d, 0, 0])^{\rm exp}$.
Since, by Fubini Theorem,
\begin{equation}\mu_{\Lambda} (p_1^* (f) p_2^* (g)) =
\mu_{\Lambda}(\mu_{x + y}(p_1^* (f) p_2^* (g))),
\end{equation}
it follows that $f \ast g$
lies in
$\cI_{\Lambda} (\Lambda [d, 0, 0])^{\rm exp}$.
\end{proof}

\begin{prop}\label{convprop}
The convolution product
$(f, g) \mapsto f \ast g$ is $\cC (\Lambda)^{\rm exp}$-linear and it endows
$\cI_{\Lambda} (\Lambda [d, 0, 0])^{\rm exp}$
with an associative and commutative law.
\end{prop}

\begin{proof}$\cC (\Lambda)^{\rm exp}$-linearity and commutativity being  clear,
let us check associativity. This follows from the fact that, if
$f$, $g$ and $h$ are be functions in
$\cI_{\Lambda} (\Lambda [d, 0, 0])^{\rm exp}$,
\begin{equation}
(f \ast g)  \ast h =  \mu_{x + y + z}(p_1^* (f) p_2^* (g) p_3^*
(h)).
\end{equation}
by Fubini  Theorem.
\end{proof}

\begin{prop}\label{convfou}Let $f$ and $g$ be two functions in
$\cI_{\Lambda} (\Lambda [d, 0, 0])^{\rm exp}$.
Then
\begin{equation*}
\Fv (f \ast g) = \Fv (f) \, \Fv (g).
\end{equation*}
\end{prop}

\begin{proof}The proof is just the same as the usual one.
Let us  consider the function $E (x (y + z)) f (y) g (z)$ on
$\Lambda [3d, 0, 0]$. It is integrable rel.~$x$, and by Fubini Theorem we
have
\begin{equation}
\begin{split}
\mu_x (E (x (y + z)) f (y) g (z)) &= \mu_x ((E (xy ) f(y)) (E (xz) g(z)))\\
&=  \Fv (f) \Fv (g).
\end{split}
\end{equation}
On the other hand, by the change of
variables formula, using the change of variables
$(x, y, z) \mapsto (x, u = y + z, z)$,
$\mu_x (E (x (y + z)) f (y) g (z))$ may be expressed as
\begin{equation}
\mu_x (E (xu) \mu_{x, u} (f (u - z) g (z))) = \Fv (f \ast g),
\end{equation}
which ends the proof.
\end{proof}

For $\varphi$ in
$\cC (\Lambda [d, 0, 0])^{\rm exp}$ we write
$\check{\varphi}$ for $\iota^* (\varphi)$, with $\iota : \Lambda [d, 0, 0]
\rightarrow \Lambda [d, 0, 0]$ the $\Lambda$-morphism sending
$x$ to $- x$.

We have the following partial Fourier inversion:

\begin{prop}\label{parfouinv}
Let $\varphi$ be a function in $\cI_{\Lambda} (\Lambda [d, 0, 0])^{\rm exp}$.
For every $\alpha$ in $\ZZ$,
$\varphi_{\alpha} \Fv (\varphi)$ lies in $\cI_{\Lambda} (\Lambda [d, 0, 0])^{\rm exp}$
and
\begin{equation*}\label{trinv}
\Fv (\varphi_{\alpha} \Fv (\varphi))=
\LL^{- \alpha d} \, \check{\varphi} \ast  \varphi_{- \alpha + 1}.
\end{equation*}
\end{prop}

\begin{proof}
We shall work on $\Lambda [3d, 0, 0]$, keeping notations and
conventions from the proof of Proposition \ref{convfou}. The
integrability of $\varphi_{\alpha} \Fv (\varphi)$ follows from the
fact that the function $E (y z) \varphi_{\alpha} (y) \varphi (z)$
on $\Lambda [2d, 0, 0]$ lies in $\cI_{\Lambda} (\Lambda [2d, 0,
0])^{\rm exp}$. Let us  consider the function $E (y (x + z))
\varphi_{\alpha} (y) \varphi (z)$ on $\Lambda [3d, 0, 0]$. It is
integrable rel.~$x$, and by Fubini Theorem we have
\begin{equation}
\begin{split}
\mu_x (E (y (x + z)) \varphi_{\alpha} (y) \varphi (z)) & = \mu_x (\mu_{(x, y)}(E (y (x + z)) \varphi_{\alpha} (y) \varphi (z)))  \\
&= \mu_x (\varphi_{\alpha}(x) E (xy) \mu_{(x, y)} (E (yz) \varphi (z)))\\
&= \mu_x (E (xy) \varphi_{\alpha}(y) \Fv (\varphi) (y))\\
&= \Fv (\varphi_{\alpha}\Fv (\varphi)).
\end{split}
\end{equation}
On the other hand, performing the change of
variables
$(x, y, z )\mapsto (u = x + z, y, z)$, we have, by the change of variables formula,
\begin{equation}
\begin{split}
\mu_x (E (y (x + z)) \varphi_{\alpha} (y) \varphi (z)) & = \mu_x (\mu_{(x, z)}(E (y (x + z)) \varphi_{\alpha} (y) \varphi (z)))  \\
&= \mu_{u - z} (\mu_{(u, z)}(E (u y) \varphi_{\alpha} (y) \varphi (z)))\\
&= \mu_{u - z} (\varphi (z)\mu_{(u, z)}(E (u y) \varphi_{\alpha} (y)))\\
&= \mu_{u - z} (\varphi (z) \Fv (\varphi_{\alpha}) (u))\\
&= \mu_{u +z} (\check{\varphi} (z) \Fv (\varphi_{\alpha}) (u))\\
&= \check{\varphi} \ast \Fv (\varphi_{\alpha}),
\end{split}
\end{equation}
which concludes the proof.
\end{proof}

\subsection{Schwartz-Bruhat functions}

We define the space
$\cS_{\Lambda} (\Lambda [d, 0, 0])^{\rm exp}$ of
Schwartz-Bruhat functions over $\Lambda$
as the $\cC (\Lambda)^{\rm exp}$-submodule
of $\cI_{\Lambda}  (\Lambda [d, 0, 0])^{\rm exp}$ consisting of functions $f$
such that
\begin{equation}\label{first}
f \cdot \varphi_{\alpha} = f  \quad \rm{for}\, \alpha \ll 0
\end{equation}
and
\begin{equation}\label{second}
f \ast \varphi_{\alpha} = \LL^{- \alpha d}\, f  \quad \rm{for}\,
\alpha \gg 0.
\end{equation}
Condition (\ref{first}) stands for ``compactly supported"
and condition (\ref{second}) for ``locally constant".
Here the quantifier  $\alpha \ll 0$ in (\ref{first}), resp.
$\alpha \gg 0$ in (\ref{second}), means there exists a definable function
$\alpha_0 : \Lambda \rightarrow \ZZ$ such that
(\ref{first}), resp.  (\ref{second}), holds for every definable function
$\alpha : \Lambda \rightarrow \ZZ$ such that $\alpha \leq \alpha_0$, resp.
$\alpha \geq \alpha_0$.

\begin{theorem}\label{sinv}Fourier transformation induces an isomorphism
\begin{equation}
\Fv : \cS_{\Lambda} (\Lambda [d, 0, 0])^{\rm exp} \simeq \cS_{\Lambda} (\Lambda [d, 0, 0])^{\rm exp}
\end{equation}
and, for every $\varphi$ in $\cS_{\Lambda} (\Lambda [d, 0, 0])^{\rm exp}$,
\begin{equation}\label{fi}
(\Fv \circ \Fv) (\varphi)
= \LL^{-d}\,  \check{\varphi}.
\end{equation}
\end{theorem}

\begin{proof}Let $\varphi$ be
in
$\cS_{\Lambda} (\Lambda [d, 0, 0])^{\rm exp}$.
Let us note that, for $\alpha \ll 0$,
\begin{equation}
\Fv (\varphi) \varphi_{\alpha} = \Fv (\varphi).
\end{equation}
Indeed,
by Proposition \ref{convfou} and Proposition \ref{form} (1),
for $\alpha \gg 0$,
\begin{equation}
\Fv (\varphi)
= \Fv (\varphi \ast \varphi_{\alpha}) \, \LL^{\alpha d}
= \Fv (\varphi) \varphi_{- \alpha + 1}.
\end{equation}
It follows from Proposition \ref{parfouinv}
that $\Fv (\varphi)$ lies in $\cI_{\Lambda} (\Lambda [d, 0, 0])^{\rm exp}$
and that
\begin{equation}
\Fv (\Fv (\varphi))=
\LL^{- \alpha d} \, \check{\varphi} \ast  \varphi_{- \alpha + 1},
\end{equation}
for $\alpha \ll 0$.
Since $\varphi$ lies in
$\cS_{\Lambda} (\Lambda [d, 0, 0])^{\rm exp}$,
$ \check{\varphi}$ also,
hence
\begin{equation}
\check{\varphi} \ast  \varphi_{- \alpha + 1}
= \LL^{(\alpha - 1)d} \, \check{\varphi}
\end{equation}
for $\alpha \ll 0$,
and we deduce (\ref{fi}).
So we are left with proving that
$\Fv (\varphi)$ lies in
$\cS_{\Lambda} (\Lambda [d, 0, 0])^{\rm exp}$.
It is enough to check that, for $\alpha \ll 0$,
\begin{equation}
\Fv (\varphi) \, \LL^{- (-\alpha + 1)d}
= \Fv (\varphi) \ast \varphi_{-\alpha +1},
\end{equation}
which follows from
the relations
\begin{equation}
\Fv (\Fv (\Fv (\varphi))) = \Fv (\varphi_{\alpha}\Fv (\Fv (\varphi)))
= \LL^{- \alpha d} \, \check{\Fv (\varphi)} \ast  \varphi_{- \alpha + 1},
\end{equation}
for $\alpha \ll 0$
by Proposition \ref{parfouinv},
 and,
by (\ref{fi}),
\begin{equation}
\Fv (\Fv (\Fv (\varphi))) = \LL^{-d} \, \check{{\Fv} (\varphi)}. 
\qedhere
\end{equation}
\end{proof}

Now we can prove Fourier inversion for integrable functions with integrable Fourier transform.

\begin{theorem}Let $\varphi$   be in
$\cI_{\Lambda} (\Lambda [d, 0, 0])^{\rm exp}$. We assume  $\Fv (\varphi)$
belongs also to $\cI_{\Lambda} (\Lambda [d, 0, 0])^{\rm exp}$. Then the functions
$(\Fv \circ \Fv) (\varphi)$ and  $\LL^{-d}\,  \check{\varphi}$ have the same class in
$C^d (\Lambda [d, 0, 0] \rightarrow \Lambda)$.
\end{theorem}

\begin{proof}By induction on $d$ we may assume $d = 1$.
Take $\varphi$ in $\cI_{\Lambda} (\Lambda [1, 0, 0])^{\rm exp}$.
By Lemma \ref{locsb} and additivity,
we may assume there exists a $1$-cell
$\lambda : Z \rightarrow Z_C \subset \Lambda [1, s, r]$,
such that, denoting by
$i$ the inclusion
$Z \rightarrow \Lambda [1, 0, 0]$ and by
$j$ the inclusion
$Z_C \rightarrow \Lambda [1, s, r]$,
$\varphi = i_! (i^* (\varphi))$
and
$\psi := j_! \lambda_! (i^* (\varphi))$
lies in
 $\cS_{\Lambda[0, r, s]}
  (\Lambda [1, r, s])^{\rm exp}$.
Since, denoting by $\pi$ the projection $\Lambda [1, r, s]
\rightarrow \Lambda$, $\varphi = \pi_! (\psi)$, the result follows
formally from Lemma \ref{pf} and Theorem \ref{sinv}.
\end{proof}

\begin{lem}\label{locsb}For every $\varphi$ in $\cI_{\Lambda} (\Lambda [1, 0, 0])^{\rm exp}$
there exists a cell decomposition of $\Lambda [1, 0, 0]$ such
that, for every $1$-cell $\lambda : Z \rightarrow Z_C \subset
\Lambda [1, s, r]$, the function $ j_!  \lambda_{!}(i^*\varphi)$
lies in $\cS_{\Lambda[0, r, s]} (\Lambda [1, r, s])^{\rm exp}$,
where $i$ denotes the inclusion $Z \rightarrow \Lambda [1, 0, 0]$
and $j$ the inclusion $Z_C \rightarrow \Lambda [1, s, r]$.
\end{lem}

\begin{proof}Follows easily from \ref{4.1},
or even from Theorem \ref{np}.
\end{proof}

\begin{lem}\label{pf}For $r$ and $s$ in $\NN$, denote by
$\pi$ the projection $\Lambda [d, r, s]\rightarrow \Lambda [d, 0,
0]$ and recall notation from \ref{kproje}. For any $\varphi$ in
$\cI_{\Lambda[0, r, s]} (\Lambda [d, r, s])^{\rm exp}$,
if the function $\pi_! (\varphi)$ lies in $\cI_{\Lambda} (\Lambda
[d, 0, 0])^{\rm exp}$, then
\begin{equation*}
\Fv (\pi_! (\varphi)) = \pi_! (\Fv (\varphi)).
\end{equation*}
\end{lem}

\begin{proof}Follows from the fact that $\mu_{p_1}$ commutes with $\pi_!$.
\end{proof}

\section{Exponential integrals over the $p$-adics}\label{sec:padics}

\subsection{Definable sets over the $p$-adics}\label{defpadic} Let $K$ be a finite field extension of $\QQ_p$ with valuation ring
$R$. We recall the notion of (globally) subanalytic subsets of $K^n$
and of semialgebraic subsets of $K^n$. Let $\cL_{\rm
Mac}=\{0,+,-,\cdot,\{P_n\}_{n>0}\}$ be the language of Macintyre and
$\cL_{\rm an}=\cL_{\rm Mac}\cup
\{^{-1},\cup_{m>0}K\{x_1,\ldots,x_m\}\}$, where $P_n$ stands for the
set of $n$th powers in $K^\times$, where $K\{x_1,\ldots,x_m\}$ is
the ring of restricted power series over $K$ (that is, formal power
series converging on
$R^m$),
 and each element $f$ of $K\{x_1,\ldots,x_m\}$ is interpreted
as the restricted analytic function $K^m\to K$ given by
\begin{equation}
x\mapsto
\begin{cases}
 f(x) & \mbox{if }x\in
R^m \\
0 & \mbox{else.}
\end{cases}\end{equation}
By subanalytic we mean $\cL_{\rm an}$-definable with coefficients
from $K$ and by semialgebraic we mean $\cL_{\rm Mac}$-definable
with coefficients from $K$. Note that subanalytic,
resp.~semialgebraic, sets can be given by a quantifier free formula
with coefficients from $K$ in the language $\cL_{\rm Mac}$,
resp.~$\cL_{\rm an}$.

In this section we let $\cL$ be either the language $\cL_{\rm
Mac}$ or $\cL_{\rm an}$ and by $\cL$-definable we will mean
semialgebraic, resp.~subanalytic when $\cL$ is $\cL_{\rm Mac}$,
resp.~$\cL_{\rm an}$. Everything in this section will hold for
both languages and we will give the appropriate references for
both languages where needed.

For each definable set $X\subset K^n$ let $\cC(X)$ be the
$\QQ$-algebra of functions on $X$ generated by functions $|f|$ and
$\ord(f)$ for all definable functions $f:X\to
K^\times$.\footnote{Instead of taking the $\QQ$-algebra we could as
well take the
$\ZZ[1/q,\{1/(1-q^i)\}_{i<0}]$-algebra, with $q$ the residue
cardinality of $K$.
}\label{foot:Qalg}
\par
For a $\cL$-definable set $X$, let $\cC^{\leq d}(X)$ be the ideal
of $\cC(X)$ generated by the characteristic functions $\11_Z$ of
$\cL$-definable subsets $Z\subset X$ of dimension $\leq d$. (For
the definition of the dimension of $\cL$-definable sets, see
\cite{sd} and \cite{DvdD}.) Note that the support of a function in
$\cC(X)$ is in general not $\cL$-definable, cf.~the function
$(x,y)\mapsto |x|-\ord(y)$ on $K\times K^\times$.
\par
By $C^d(X)$ we denote the quotient
\begin{equation}
C^d(X) := \cC^{\leq d}(X)  / \cC^{\leq d-1}(X).
\end{equation}
Finally we set
\begin{equation}
C(X):=\bigoplus_{d\geq 0}C^d(X).
\end{equation}
It is a module over $\cC(X)$. If $\varphi$ is in $\cC(X)$ with
support contained in a $\cL$-definable set of dimension $d$, we denote by $[\varphi]_d$ its image
in $C^d(X)$.
\par

\subsection{The $p$-adic measure}\label{pmes}
Suppose that $X\subset K^n$ is a $\cL$-definable set of dimension
$d\geq 0$. The set $X$ contains a definable nonempty  open
submanifold $X'\subset K^n$ such that $X\setminus X'$ has
dimension $<d$, cf.~\cite{DvdD}. There is a canonical
$d$-dimensional measure on $X'$ coming from the embedding in
$K^n$, which is constructed as follows, cf. \cite{serre}. For each
$d$-element subset $J$ of $\{1,\ldots,n\}$, with $j_i<j_{i+1}$,
$j_i$ in $J$, let $dx_J$ be the $d$-form
$dx_{j_1}\wedge\ldots\wedge dx_{j_d}$ on $K^n$, with
$x=(x_1,\ldots,x_n)$ standard global coordinates on $K^n$. Let
$x_0$ be a point on $X'$ such that $x_I$ are local coordinates
around $x_0$ for some $I\subset\{1,\ldots,n\}$. For each
$d$-element subset $J$ of $\{1,\ldots,n\}$ let $g_J$ be the
$\cL$-definable function determined at a neigborhood of $x_0$ in
$X'$ by $g_Jdx_I=dx_J$. There is a unique volume form
$\vert \omega_{0}\vert_{X'}$ on $X'$ which is  locally equal to
$(\max_J|g_J|) \vert dx_I \vert$ around every point $x_0$ in $X'$.
Indeed, $\vert \omega_{0}\vert_{X'}$ is equal to $\sup_{J}\vert dx_J
\vert$. The canonical $d$-dimensional measure on $X'$, cf.
\cite{serre} \cite{oesterle}, is the one induced by the volume form
$\vert \omega_{0}\vert_{X'}$. We extend this measure to $X$ by zero
and denote it by $\mu^d$.
\par
This measure allows us to define the subgroup ${\rm I}C^d(X)$ of
$C^d(X)$ for a $\cL$-definable set $X$ of dimension $d$, as the
group generated by elements $[\varphi]_d$ with $\varphi$ in $\cC(X)$
integrable for $\mu^d$. We define ${\rm I} C^e(X)$ for general $e$
as the subgroup of $C^e(X)$ consisting of elements $[\varphi]_e$
with $\varphi$ with support contained in a $\cL$-definable subset
$Z\subset X$ of dimension $e$ and with $[\varphi_{|Z}]_e$ in ${\rm
I} C^e(Z)$. Finally, we define the graded group ${\rm I} C(X)$ as
$\oplus_r {\rm I} C^r(X)$.

\subsection{Jacobian}\label{jac}Using the pullback of differential forms under analytic maps, it
is possible to define the norm of the Jacobian $|\Jac f|$ of a $\cL$-definable bijection
$f:X\subset K^n\to Y\subset K^m$ as follows.
There exists definable $K$-analytic
manifolds $X'\subset X$ and $Y'\subset Y$ such that $X\setminus
X'$ and $Y\setminus Y'$ have dimension
$<d$ with $d=\dim X$
 and such that $f_{|X'}$ is a $K$-bi-analytic bijection
onto $Y'$. For subsets $I$ and $J$ of $\{1, \dots, n\}$ and $\{1,
\dots, m\}$ respectively, we denote by $U_{I, J}$ the definable
subset of $X'$ consisting of points $x_0$  such that $\vert dx_I
\vert$ coincides with
 $\vert \omega_0 \vert_{X'}$
on a
neighborhood of $x_0$ and
$\vert dy_J \vert$ coincides with
 $\vert \omega_0 \vert_{Y'}$
on a
neighborhood of $f(x_0)$.
On $U_{I, J}$ we may write
$f^* (dy_J) = g_{I, J} dx_I$.
The functions $\vert g_{I, J} \vert$ are constructible on $U_{I, J}$
and there exists a unique constructible function $h$ on $X'$
restricting to $\vert g_{I, J} \vert$ on each $U_{I, J}$. We define
$|\Jac f|$ as the class of $h$ in $C^d(X)$ which does not depend on
the choices made.

The $p$-adic change of variables formula (cf. \cite{Igusa}) may 
be
restated as follows:

\begin{prop}Let
$f:X\subset K^n\to Y\subset K^m$ be a $\cL$-definable bijection,
with $d=\dim X$.
 For every measurable subset $A$ of $Y$ one has
\begin{equation}
\mu^d(A)= \int_{f^{-1}(A)} |\Jac f|\,  \mu^d.
\end{equation}
\end{prop}

\subsubsection{}\label{delta} For the proof of Theorem \ref{mainK} below
we shall need the following variant of $|\Jac f|$.
Let $X$ a definable subset of $Y \times K^n$ for $Y$ a definable
subset of $K^m$ and consider the morphism $f : X \rightarrow Y$
induced by projection. Assume first $X$ is of dimension $r$, $Y$ of
dimension $s$, $f$ is surjective and all fibers of $f$ have
dimension $r - s$. In this setting we define a constructible
function $\delta$ defined almost everywhere on $X$ as follows. We
can choose (cf.~\cite{DvdD})
definable  
manifolds $X'\subset X$ and $Y' \subset Y$ such that $X\setminus X'$
has dimension $<r$, $Y\setminus Y'$ has dimension $<s$, $f$
restricts to a locally analytic morphism $f' : X' \rightarrow Y'$
and for every point $y$ in $Y'$,
$f^{-1} (y) \setminus f'{}^{-1}(y)$ is of dimension $< r - s$, and
such that $f'$ is regular (that is, the Jacobian matrix has
everywhere maximal rank).
 Denote by $y_i$ the coordinates on $K^m$ and by $z_i$ the
coordinates on $K^n$. Consider subsets $I$, $I'$ of $\{1, \dots,
m\}$ and $J$, $J'$ of $\{1, \dots, n\}$, 
and denote by $U_{I, J, I', J'}$
 the set of points $x$ of $X'$
such that on a neighborhood of $x$, $\vert \omega_0 \vert_{X'}$
coincides with $\vert dy_I \wedge dz_J \vert$, $\vert \omega_0
\vert_{f'{}^{-1}(f(x))}$ coincides with $\vert dz_{J'}\vert$, and on
a neighborhood of $f (x)$, $\vert \omega_0 \vert_{Y'}$ coincides
with $\vert dy_{I'}\vert$. On  $U_{I, J, I', J'}$ we may  write
\begin{equation}
dy_{I'} \wedge dz_{J'}
= g_{I, J, I', J'}
dy_I \wedge dz_J,
\end{equation}
with $g_{I, J, I', J'}$ definable.
There is a unique constructible function $g$ on $X'$
restricting to $\vert g_{I, J, I', J'} \vert$ on $X'$. We denote its class in
$C^r (X)$ , which is independent of the choices made, by $\delta (f)$.
Note that when $f$ is an isomorphism $\delta (f) = \vert \Jac f\vert$.

The proof of the following chain rule is clear:

\begin{lem}\label{trans}Let $Z$ be a definable subset of dimension $t$ of $K^m$,
$Y$ be a definable subset of dimension $s$ of $Z \times K^n$ and $X$
be a definable subset of dimension $r$ of $Y \times K^q$ Assume the
morphisms $f : X \rightarrow Y$ and $g : Y \rightarrow Z$ are
induced by projections and are surjective, and that all their fibers
have dimension $r - s$ and $s- t$, respectively. Then the equality
$$
\delta (g \circ f) = \delta (f) (\delta (g) \circ f)
$$
holds. \qed
\end{lem}

\subsection{$p$-adic Cell Decomposition}
Cells are defined by induction on the number of variables:
 \begin{definition}\label{def::cell}
A $\cL$-cell $A\subset K$ is a (nonempty) set of the form
 \begin{equation}
 \{t\in K\mid |\alpha|\sq_1 |t-c|\sq_2 |\beta|,\
  t-c\in \lambda P_n\},
\end{equation}
with constants $n>0$, $\lambda,c$ in $ K$, $\alpha,\beta$ in
$K^\times$, and $\square_i$ either $<$ or no condition. A $\cL$-cell
$A\subset K^{m+1}$, $m\geq0$, is a set  of the form
 \begin{equation}\label{Eq:cell:decay}
 \begin{array}{ll}
\{(x,t)\in K^{m+1}\mid
 &
 x\in D, \  |\alpha(x)|\sq_1 |t-c(
 x)|\sq_2 |\beta(x)|,\\
 &
  t-c(x)\in \lambda P_n\},
 \end{array}
 \end{equation}
 with $(x,t)=(x_1,\ldots,
x_m,t)$, $n>0$, $\lambda$ in $ K$, $D=\pi_m(A)$ a cell where $\pi_m$
is the projection $K^{m+1}\to K^m$, $\cL$-definable functions
$\alpha,\beta:K^m\to K^\times$ and $c:K^m\to K$, and $\square_i$
either $<$ or no condition, such that the functions
$\alpha,\beta$, and $c$ are analytic on $D$. We call $c$ the
center of the cell $A$ and $\lambda P_n$ the coset of $A$. In
either case, if $\lambda=0$ we call $A$ a $0$-cell and if
$\lambda\not=0$ we call $A$ a $1$-cell.
 \end{definition}

In the
$p$-adic semialgebraic case, Cell Decomposition Theorems are  due to Cohen \cite{cohen}
and Denef \cite{D84}, \cite{Dcell} and they were extended in
\cite{Ccell}
to the subanalytic setting where one can find the
following version:

\begin{theorem}[$p$-adic Cell Decomposition]\label{thm:CellDecomp}
Let $X\subset K^{m+1}$ and $f_j:X\to K$ be $\cL$-definable for
$j=1,\ldots,r$. Then there exists a finite partition of $X$ into
$\cL$-cells $A_i$ with center $c_i$ and coset $\lambda_i P_{n_i}$
such that
 \begin{equation*}
 |f_j(x,t)|=
 |h_{ij}(x)|\cdot|(t-c_i(x))^{a_{ij}}\lambda_i^{-a_{ij}}|^\frac{1}{n_i},\quad
 \mbox{ for each }(x,t)\in A_i,
 \end{equation*}
with $(x,t)=(x_1,\ldots, x_m,t)$, integers $a_{ij}$, and
$h_{ij}:K^m\to K$ $\cL$-definable functions which are analytic on
$\pi_m(A_i)$, $j=1,\ldots,r$. If $\lambda_i=0$, we use the
convention that $a_{ij}=0$.
 \end{theorem}

We shall also use the following lemma from \cite{Cexp}:
\begin{lem}\label{prop:descrip:simple}
Let $X\subset K^{m+1}$ be $\cL$-definable and let $G_j$ be
functions in $\cC(X)$ in the variables $(x_1,\ldots,x_m,t)$ for
$j=1,\ldots,r$. Then there exists a finite partition of $X$ into
$\cL$-cells $A_i$ with center $c_i$ and coset $\lambda_i P_{n_i}$
such that each restriction $G_j|_{A_i}$ is a finite sum of
functions of the form
\[
|(t-c_i(x))^{a}\lambda^{-a}|^\frac{1}{n_i}v(t-c_i(x))^{s}h(x),
\]
where $h$ is in $\cC(K^m)$, and $s\geq 0$ and $a$ are integers.
\end{lem}

\subsection{Integration}

By $\Def(\cL)$ we denote the category of 
$\cL$-definable subsets $X\subset K^n$ for $n>0$, with
$\cL$-definable maps as morphisms. We can now state a general
integration result which states unicity and existence of a certain
integral operator. This integral operator is introduced as a
push-forward operator of functions under $\cL$-definable maps,
inspired by integration in the fibers with a measure on the fibers
coming from Leray-differential forms.

In \cite{D2000}, see also \cite{D84}, 
Denef proved stabillity of $p$-adic constructible
functions under integration with respect to parameters in the
semialgebraic case. Denef's result  had a major influence on our
work \cite{cl} and the present one. It has been later generalized to
the subanalytic case by the first author in  \cite{Ccell} and
\cite{Cexp}.

\begin{prop}[\cite{D2000}, \cite{Ccell}, \cite{Cexp}]\label{linear}Let $W$ be a definable subset of $K^n$ of dimension $r$.
\begin{enumerate}
\item[(1)]
Let $\varphi$ be in $\cC(W \times K^{m})$. Assume for every $x$ in
$W$ the function
$t \mapsto \varphi (x, t)$
 is integrable on $K^m$. Then the function
 \begin{equation*}
 g(x) : =\int_{K^m}  \varphi(x,t)|dt|
 \end{equation*}
 lies in $\cC (W)$.
 \item[(2)]Let $\varphi$ be in ${\rm I}C^{r+m} (W \times K^m)$. Then, there exists
a function $g$ in $\cC(K^n)$ such that for all $x$ in $W\setminus
Z$, with $Z$ a $\cL$-definable set of dimension $< r$ in  $K^n$, one
has
 \begin{equation*}
 g(x)=\int_{K^m}
  \varphi(x,t)|dt|.
   \end{equation*}
 \end{enumerate}
\end{prop}

\begin{proof}When
$W = K^n$, statement (1) is proved in
\cite{D2000} in the semialgebraic case and in
\cite{Ccell} in the subanalytic case and statement (2) is proved
in
\cite{Cexp}. The proofs carry  over litterally to general $W$.
\end{proof}

\begin{remark}
The point in (2) is that it is possible that the set of points in
$W$ such that the restriction of $\varphi$ to the fibers
over which it is
integrable may not be definable.
\end{remark}

We shall now prove the following
analogue  of Theorem 10.1.1 of \cite{cl}. Note that the proof will be much easier, since
integrable functions
are already defined and Proposition \ref{linear}
is available.

\begin{theorem}\label{mainK}
There exists a unique functor from $\Def(\cL)$ to the category of
groups sending a $\cL$-definable set $X$ to the group ${\rm I} C(X)$ such
that a morphism $f:X\to Y$ in $\Def(\cL)$ is sent to a group
morphism $f_!: {\rm I} C(X)\to {\rm I} C(Y)$ satisfying the following axioms

\medskip
\item[]{\textup{(A1) Disjoint union: }} Assume that $X$, resp.~$Y$, is the
disjoint union of two $\cL$-definable sets $X_1$ and $X_2$,
resp.~$Y_1$ and $Y_2$, such that $f(X_i)\subset Y_i$. Write
$f_i:X_i\to Y_i$ for the restrictions. Then we have
$f_! = f_{1 !} \oplus f_{2!}$
under the isomorphisms ${\rm I}C (X) \simeq {\rm I} C (X_1) \oplus
{\rm I} C (X_2)$ and ${\rm I}C (Y) \simeq {\rm I} C (Y_1) \oplus
{\rm I} C (Y_2)$.

\medskip
\item[]{\textup{(A2) Projection formula: }}For every $\alpha$ in $\cC
(Y)$ and $\beta$ in ${\rm I}C (X)$, if $(\alpha\circ f) \beta$ is in
${\rm I} C (X)$, then $f_!
((\alpha\circ f) \beta) = \alpha f_! (\beta)$.

\medskip
\item[]{\textup{(A3) Projection for $1$-cells: }} Let
$X\subset K^{n+1}$
 be a $1$-cell of dimension $r$ and $Y$ its image
under the projection on $K^n$, $f:X\to Y$ the projection. Let
$\varphi$ be a $\mu^r$- integrable function in $\cC (X)$. By
Proposition \ref{linear} there exists a $\cL$-definable set
$Z\subset Y$ such that $Y\setminus Z$ has dimension $<r-1$ and such
that the the function $g:Y\to\QQ:y\mapsto \int_{f^{-1} (y)} \11_{Y
\setminus Z}(y) \varphi(y,t)\vert dt\vert$ lies in $\cC^{\leq
r-1}(Y)$. We let $f_!([\varphi]_r)$ be the class of $g$ in ${\rm
I}C^{r-1}(Y)$.

\medskip
\item[]{\textup{(A4) Projection for $0$-cells: }} Let
$X\subset K^{n+1}$
 be a
$0$-cell of dimension $r$ and $Y$ its image under the projection
on $K^n$, $f:X\to Y$ the projection. Then, we define $f_!(\11_X)$
as $(|(\Jac f)\vert \circ f^{-1})^{-1}\11_Y$ in ${\rm I}C^{r}(Y)$, where
$\vert \Jac f \vert$ is as in \ref{jac}.
\end{theorem}

\begin{proof}
We will freely use classical forms of the change of variables
formula, without mentioning it.
Let us first check unicity. Since, by the graph construction, any
morphism $f: X \rightarrow Y$ is the composition of the graph
injection
$i_f : X \rightarrow X\times Y$
  and the  projection $p: X \times Y
\rightarrow Y$, it is enough, by functoriality, to prove unicity for
$i_{f!}$ and $p_!$. For projections
$X\times Y\to Y$,
 one can assume
$X = K^m$ and $Y = K^n$ by (A1). By induction on $m$, it is enough
to define $p_!$ when $m = 1$. Consider $\varphi$ in $\cC^r (K^{n +
1})$ and assume it is integrable. By cell decomposition and
linearity we may assume the support of $\varphi$ is contained in a
$Z$ cell of dimension $r$. If $Z$ is a $1$-cell, $f_!([\varphi]_r)$
is given by (A3). In case $Z$ is a $0$-cell, we may assume by (A2)
that $\varphi = \11_Z$, and then $f_!([\varphi]_r)$ is given by
(A4). Finally, since $q \circ i_f = {\rm id}_X$, with $q$ the
projection $X \times Y \rightarrow X$, and since $q$ induces a
bijection between the graph of $f$ and $X$, unicity for  $i_{f!}$
reduces to that of $q_!$ (an essentially similar argument is given
with full details in the uniqueness section of the proof  of Theorem
10.1.1 of \cite{cl}).

Let us now define $f_!$ for projections. Let $X$ be a definable
subset of $Y \times K^n$ for $Y$ a definable subset of $K^m$ and
consider the morphism $f : X \rightarrow Y$ induced by projection.
Assume first that $X$ is of dimension $r$, $Y$ of dimension $s$, $f$
is surjective and all fibers of $f$ have dimension $r - s$.

Let $\varphi$ be a
$\mu^r$-integrable
 function in $\cC (X)$. There exists a definable
subset $Z$ of $Y$, with dimension $< s$, such that the function
\begin{equation}\label{pmo}
g:y \mapsto \int_{f^{-1}(y)} \11_{Y \setminus Z}(y)
  \delta \varphi \mu^{r-s}
\end{equation}
lies in $\cC (Y)$,
with $\delta$ as in section \ref{delta}, and $\mu^{r-s}$ the measure
as in section \ref{pmes}.
 Indeed, by Fubini and induction, and
possibly after considering a partition of $X$ and $Y$, we may assume
$n = 1$. Then, we may by cell decomposition assume $X$ is a cell. If
$X$ is a $0$-cell the statement is clear, if $X$ is a $1$-cell the
statement follows from Proposition \ref{linear}. We may then define
$f_! ([\varphi]_r)$ to be the class of $g$ in
$C^{r-s} (Y)$.
 It follows from Fubini that $f_! ([\varphi]_r)$ lies
in ${\rm I} C^{r - s} (Y)$. Note that certainly (A3) and (A4) hold.
Also, there is a unique way to extend that construction to a
morphism $f_! : {\rm I} C (X) \rightarrow {\rm I} C (Y)$ satisfying
(A1) and (A2)
for every morphism $f : X \rightarrow Y$ induced by a projection $Y
\times K^m \rightarrow Y$. Indeed, it is enough to construct $f_!$
on ${\rm I} C^r (X)$ and after cutting $X$ and $Y$ into finitely
many definable pieces, one may assume the above condition is
verified. Furthermore, by Lemma \ref{trans}, $(g \circ f)_! = g_!
\circ f_!$ for composable morphisms induced by projections.

Let us now define
$i_!$ when $i : X \rightarrow Y$ is a definable injection. 
Let $\varphi$ be in $C^e (X)$. Consider a definable subset $X'$ of dimension $e$ of $X$ such that
the support of a representative $\tilde \varphi$ of $\varphi$
is contained in $X'$.
Denote by $\lambda : i (X') \rightarrow X'$ the inverse of the isomorphism induced by $i$.
We define
$i_! (\varphi)$ as the image
of  $[(\tilde \varphi \circ \lambda )]_e \vert \Jac (\lambda)\vert$
in $C^e (Y)$ under the inclusion  $C^e (i (X')) \rightarrow C^e (Y)$.
Certainly $i_! (\varphi)$ is integrable if $\varphi$ is, hence
we deduce a morphism
$\lambda_! : {\rm I} C (X) \rightarrow {\rm I} C (Y)$.

For a general morphism $f : X \rightarrow Y$ one considers the
factorisation $f = \pi_f \circ i_f$, with $i_f : X \rightarrow X
\times Y$ the inclusion of the graph and $\pi_f : X \times Y
\rightarrow Y$ the projection, and one sets $f_! := \pi_{f !}\circ
i_{f !}$. One is then left with checking that the construction
coincides with the previous one for injections and projections and
that $(g \circ f)_! = g_! \circ f_!$ for composable morphisms. This
is purely formal and performed exactly as the proof of the
corresponding statements in the proof of Theorem 10.1.1 and
Proposition 12.1.2 of \cite{cl}.
\end{proof}

\subsection{Exponential constructible functions}

Fix an additive character $\psi:K\to\CC^\times$ which is trivial on
the maximal ideal $M$ of $R$ and such that $\psi(x)\not=1$ for some
$x$ in $K$ with $\ord(x)=0$.

For $X$ a $\cL$-definable set, we let $\cC(X)^{\rm exp}$ be the
$\QQ$-algebra generated by $\cC(X)$ and all functions $\psi(f)$,
where $f:X\to K$ is $\cL$-definable (cf.~footnote \ref{foot:Qalg}).

Similarly, for each $d\geq 0$ we define $\cC ^{\leq d} (X)^{\rm exp}
$ as the $\QQ$-algebra generated by $\cC^{\leq d}(X)$ and all
functions $\psi(f)$ with $f:X\to K$ $\cL$-definable.

We set
\begin{equation}C (X)^{\rm exp} = \oplus_d  C^d  (X)^{\rm exp}\end{equation}
 with
\begin{equation}C^d  (X)^{\rm exp} := \cC ^{\leq d} (X)^{\rm exp} / \cC ^{\leq
d-1} (X)^{\rm exp}.
 \end{equation}
 We call elements of $C (Z)^{\rm exp}$
constructible exponential Functions.

For $d\geq0$ we define the group $ {\rm I} C^d (X)^{\rm exp}$ of
integrable constructible exponential Functions as the
$\QQ$-algebra of functions which are defined almost everywhere,
generated by ${\rm I} C (X)$ and all functions $\psi(f)$ with
$f:X\to K$ $\cL$-definable. It is a graded subgroup of $C (Z)^{\rm
exp}$ and a module over $\cC (Z)^{\rm exp}$.

The following exponential analogue of Proposition \ref{linear} is
our main $p$-adic result:

\begin{prop}\label{linearexp}
Let $W$ be a definable subset of $K^n$ of dimension $r$.
\begin{enumerate}
\item[(1)]
Let $\varphi$ be in $\cC(W \times K^{m})^{\rm exp}$. Assume for
every $x$ in $W$ the function $t \mapsto \varphi (x, t)$ is
integrable on $K^m$. Then the function
 \begin{equation*}
 g(x) : =\int_{K^m}  \varphi(x,t)|dt|
 \end{equation*}
 lies in $\cC (W)^{\rm exp}$.
 \item[(2)]Let $\varphi$ be in ${\rm I}C^{r+m} (W \times K^m)^{\rm exp}$. Then, there exists
a function $g$ in $\cC(K^n)^{\rm exp}$ such that for all $x$ in
$W\setminus Z$,
with $Z$ a $\cL$-definable set of dimension $< r$ in $K^n$, one has
 \begin{equation*}
 g(x)=\int_{K^m}
  \varphi(x,t)|dt|.
   \end{equation*}
 \end{enumerate}
\end{prop}

\begin{proof}(1) follows easily from (2), so let us prove (2).
By Fubini Theorem it is enough to consider the case $m=1$. By
linearity of the integral operator it is enough to prove the
Proposition when  $\varphi=\varphi_0 \psi(f)$ with $\varphi_0$ in
${\rm I}C^{r+1}(W \times K)$ and $f:W \times K\to K$ a $\cL$-definable morphism.

We partition $W \times K$ into $\cL$-definable parts $B_1$, $B_2$, and
$B_3$:
\begin{eqnarray*}
 B_1 & := &
 \{(x,t)\in W \times K\mid f(x,\cdot) \mbox{ is }C^1
 \mbox{ at $t$ and }\frac{\partial f}{\partial t}(x,t)\not=0 \},
 \\
 B_2 & := & \{(x,t)\in W \times K\mid f(x,\cdot) \mbox{ is not }C^1
 \mbox{ at $t$}\},
 \\
 B_3 & := & \{(x,t)\in W \times K\mid f(x,\cdot) \mbox{ is }C^1
\mbox{ at $t$ and }\frac{\partial f}{\partial t}(x,t)=0 \},
\end{eqnarray*}
where $C^1$ at a point means continuously differentiable in an
open neighborhood and $f(x,\cdot)$ denotes the function $K\to
K:t\mapsto f(x,t)$ for each $x$ in $W$.

We can partition $B_2$ into finitely $\cL$-definable sets $B_{2i}$
such that $f(x,t)=g_i(x)$ for each $i$ and for each $(x,t)$ in
$B_{2i}$ for some $\cL$-definable functions $g_i:W\to K$. This
follows from the fact  that for each $x$ in $W$ the set
$B_{2x}:=\{t\in K\mid (x,t)\in B_{2}\}$ is finite and uniformly
bounded in number when $x$ varies, and from the $\cL$-definability
of choice functions, that is, there exist $\cL$-definable sections
of $\cL$-definable maps (these statements follow immediately from
cell decomposition \ref{thm:CellDecomp}, cf.~also \cite{DvdD} and
\cite{sd}).

Similarly, we can partition $B_3$ into finitely many
$\cL$-definable sets $B_{3i}$ such that $f(x,t)=r_i(x)$ for each
$i$ and for each $(x,t)$ in $B_{3i}$ for some $\cL$-definable
functions $r_i:W\to K$. This follows from the facts that for
each $x$ the image of the function $t\mapsto f(x,t)$ is discrete,
hence finite and uniformly bounded when $x$ varies, and from the
$\cL$-definability of choice functions (again these statements
follow from cell decomposition \ref{thm:CellDecomp}, cf.~also
\cite{DvdD} and \cite{sd}).

Hence, for the functions $\11_{B_{\ell i}}\varphi$ with $\ell=2,3$
the Proposition follows. By linearity of the integral operator we
only have to prove the Proposition for the function
$\11_{B_1}\varphi$.

By the implicit function theorem, the set $\{t\mid f(x,t)=z,\
(x,t) \in B_1\}$ is discrete for each $x$ in $W$ and $z$ in $K$,
hence finite and uniformly bounded when $x$ and $z$ vary, by the
cell decomposition \ref{thm:CellDecomp} (or by \cite{DvdD},
\cite{sd}). By the existence of $\cL$-definable choice functions,
we can partition $B_1$ into finitely many $\cL$-definable parts
$B_{1i}$ such that $f(x,\cdot)$ is injective on $B_{1ix}$ for each
$i$ and each $x$, with $B_{1ix}:=\{t\in K\mid (x,t)\in B_{1i}\}$.
Hence, we may as well suppose that $f(x,\cdot)$ is injective on
$B_{1x}:=\{t\in K\mid (x,t)\in B_{1}\}$ itself. Then, we let $T$
be the transformation
\begin{equation}
T:
\begin{cases}
B_1\mapsto T(B_1) &\\
(x,t)\mapsto (x,y):=(x,f(x,t)),&
\end{cases}
\end{equation}
and let $|\Jac T|$ be the Jacobian of $T$ as in \ref{jac}. Writing
$T(B_1)_x$ for $\{t\in K\mid (x,t)\in T(B_{1})\}$, one has by the
change of variables rule for each $x$ in $ \pi_{B_1}(B_1)$
\begin{equation}
\int_{B_{1x}}\varphi_0(x,t)\psi(f(x,t))|dt| =
\int_{T(B_1)_x}(|\Jac T|\circ
T^{-1}(x,y))^{-1}\varphi_0(T^{-1}(x,y))\psi(y)|dy|.
\end{equation}
Now apply Lemma \ref{prop:descrip:simple} to the function
\begin{equation}
 \varphi_1:
  \begin{cases}T(B_1)\to \QQ &\\
  (x,y)\mapsto(|\Jac T|\circ
T^{-1}(x,y))^{-1}\varphi_0( T^{-1}(x,y))&
\end{cases}
 \end{equation}
with respect to the variable $y$ to obtain a partition of $T(B_1)$
into $\cL$-cells $A$ with center $c$ and coset $\lambda P_m$ such
that each $\varphi_{1|A}$ is a finite sum of functions of the form
\begin{equation}\label{eq:term:H}
 H(x,y)=|(y-c(x))^a\lambda^{-a}|^\frac{1}{m}v(y-c(x))^{s}h(x),
\end{equation}
where $h: W \to \QQ$ is in $\cC(W)$, and $s\geq0$ and $a$ are
integers.

\begin{claim}Possibly after refining the partition, we
can assure that for each $A$ either the projection
$A':=\pi_{W}(A)\subset W$ has zero
$\mu^r$-measure,
or we can write $\varphi_{1|A}$ as a sum of terms $H$ of the form
(\ref{eq:term:H}) such that $H$ is
$\mu^r$-integrable
over $A$ and $H(x,\cdot)$ is integrable over $A_{x}:=\{y\mid
(x,y)\in A\}$ for all $x$ in $A'$.
\end{claim}

As this claim is very similar to Claim 2 of \cite{Cexp} we will only
give an indication of its proof. Partitioning further, we may
suppose that $v(y-c(x))$ either takes only one value on $A$ or
infinitely many values. If $v(y-c(x))$ only takes one value on $A$,
we may suppose that the exponents $a$ and $s$ as in
(\ref{eq:term:H}) are zero. Now apply Lemma
\ref{prop:descrip:simple} to each $h$ and to the norms of all the
$\cL$-definable functions appearing in the description of the cells
$A$ in a similar way and do this inductively for each variable. This
way, the claim is reduced to a summation problem over a Presburger
set of integers, which is easily solved. This proves the claim.

Fix a cell $A$ and a term $H$ as in the claim. The cell $A$ has by
definition the following form
\[
\begin{array}{ll}
A=\{(x,y)\mid & x\in A',\ v(\alpha(x))\sq_1 v(y-c(x))\sq_2
v(\beta(x)),
\\
 &
 y-c(x)\in\lambda P_m \},
\end{array}
\]
where $A'=\pi_{W}(A)$ is a cell, $\sq_i$ is $<$ or no condition, and
$\alpha,\beta: W \to K^\times$ and $c: W \to K$ are $\cL$-definable
functions. We focus on a cell $A$
of dimension $r+1$,
in particular, $\lambda\not=0$ and $A'$ is of dimension $r$.

For $x$ in $A'$, we denote by $I(x)$ the value
 \[
I(x)=\int_{y\in A_{x}} H(x,y)\, \psi(y)\,|dy|,
 \]
where $A_x=\{y\in K\mid (x,y)\in A\}$. Write
 \begin{equation}G(j):=\int_{v(u)=j,\ u\in \lambda P_m }\psi(u)\,|du|.
 \end{equation}
We easily find
\begin{equation}\label{eq:sum0}
\begin{array}{c}
I(x)=\psi(c(x))\
h(x)|\lambda|^{-a/m}\sum\limits_{(\ref{summation})} q^{-ja/m}\,
j^{s} G(j),
\end{array}
\end{equation}
where the summation is over
 \begin{equation}\label{summation}
J:= \{j \mid v(\alpha(x)) \sq_1 j \sq_2 v(\beta(x)),\ j\equiv
v(\lambda)\bmod m\}.
\end{equation}
\par
By Hensel's Lemma, there exists an integer $e\geq 0$ such that all
units $u$ with $u\equiv 1\bmod \pi_0^e$ are $m$-th powers (here,
$\pi_0$ is such that $v(\pi_0)=1$). Hence, $G(j)$ is zero whenever
$j\leq-e$ (since in this case one essentially sums a nontrivial
character over a finite group). Also, when $j>0$ then
$G(j)=\int_{v(u)=j,\ u\in \lambda P_m }|du|$, which is independent
of $\psi$. We find that $I(y)$ is equal to
\begin{equation}\label{eq:sum00}
 \psi(c(y))\ h(y)|\lambda|^{-a/m}\cdot 
\big(\sum_{\sur{-e\leq j \leq 0}{j\in J}}
 q^{-ja/m}\, j^{s} G(j)
    + \sum_{\sur{0<j}{j\in J}}
    q^{-ja/m}\, j^{s}G(j)\big).
\end{equation}
The factors of (\ref{eq:sum00}) before the brackets clearly are in
$\cC(W)^{\rm exp}$. The (parameterized) finite sum inside the
brackets of (\ref{eq:sum00}) can be written as a finite sum of
generators of $\cC(W)^{\rm exp}$ since each
$G(j)=p^{-j}\alpha_{j\bmod e+n}$ with each $\alpha_{j\bmod e+n}$
some $\QQ$-linear combination of values of $\psi$ which only
depends on $j\bmod e+n$, and hence, it is also in $\cC(W)^{\rm
exp}$. The infinite sum inside the brackets of (\ref{eq:sum00}) is
in $\cC(W)$ by Proposition \ref{linear} and the above
discussion. This finishes the proof of the Proposition.
\end{proof}

One may extend Theorem \ref{mainK}
to the exponential setting as follows:
\begin{theorem}\label{mainKexp}
There exists a unique functor from $\Def(\cL)$ to the category of
groups sending a $\cL$-definable set $X$ to the group ${\rm I} C(X)^{\rm
exp}$ such that a morphism $f:X\to Y$ in $\Def(\cL)$ is sent to a
group morphism $f_!: {\rm I} C(X)^{\rm exp}\to  {\rm I} C(Y)^{\rm exp}$ satisfying
the following axioms

\medskip
\item[]{\textup{(A1) Compatibility:}} For every morphism $f : X
\rightarrow Y$ in $\Def(\cL)$, the map $f_!: {\rm I} C (X)^{\rm
exp} \rightarrow {\rm I} C  (Y)^{\rm exp}$ is compatible with the
inclusions of groups ${\rm I} C  (X) \to {\rm I} C (X)^{\rm exp}$
and ${\rm I} C  (Y) \to {\rm I} C (Y)^{\rm exp}$ and with the map
$f_!:{\rm I} C  (X)\to {\rm I} C  (Y)$ as constructed in Theorem
\ref{mainK}.

\medskip
\item[]{\textup{(A2) Disjoint union: }} Assume that $X$, resp.~$Y$, is the
disjoint union of two $\cL$-definable sets $X_1$ and $X_2$,
resp.~$Y_1$ and $Y_2$, such that $f(X_i)\subset Y_i$. Write
$f_i:X_i\to Y_i$ for the restrictions. Then we have $f_! = f_{1 !}
\oplus f_{2!}$ under the isomorphisms ${\rm I}C (X)^{\rm exp} \simeq
{\rm I} C (X_1)^{\rm exp} \oplus {\rm I} C (X_2)^{\rm exp}$ and
${\rm I}C (Y)^{\rm exp} \simeq {\rm I} C (Y_1)^{\rm exp} \oplus {\rm
I} C (Y_2)^{\rm exp}$.

\medskip
\item[]{\textup{(A3) Projection formula: }}For every $\alpha$ in $\cC
(Y)^{\rm exp}$ and $\beta$ in ${\rm I} C (X)^{\rm exp}$,
if $(\alpha\circ f) \beta$ is in ${\rm I}C (X)^{\rm exp}$,
then $f_! ((\alpha\circ f) \beta) = \alpha f_! (\beta)$.

\medskip
\item[]{\textup{(A4) Projection for $1$-cells: }} Let $X\subset K^{n+1}$ be a
$1$-cell of dimension $r$ and $Y$ its image under the projection to
$K^n$, $f:X\to Y$ the projection.  Let $\varphi$ be a $\mu^r$-
integrable function in $\cC (X)^{\rm exp}$. By Proposition
\ref{linearexp} there exists a $\cL$-definable set $Z\subset Y$ such
that $Y\setminus Z$ has dimension $<r-1$ and such that the the
function $g:Y\to\QQ:y\mapsto \int_{f^{-1} (y)} \11_{Y \setminus
Z}(y) \varphi(y,t)\vert dt\vert$ lies in $\cC^{\leq r-1}(Y)^{\rm
exp}$. We let $f_!([\varphi]_r)$ be the class of $g$ in ${\rm
I}C^{r-1}(Y)^{\rm exp}$.
\end{theorem}

\begin{proof}The proof is quite formal and similar to proofs we have already given.
Indeed, uniqueness is proven along similar lines than those  given in \ref{munique},
and for existence one can proceed similarly as in the proof of Theorem \ref{mainK},
using Proposition \ref{linearexp} instead of Proposition \ref{linear}.
\end{proof}

\subsection{Variants: adding sorts and relative versions}\label{sr}
By analogy with the motivic framework, we now expand the language
$\cL$ to a three sorted language $\cL'$ having  $\cL$ as language
for the valued field sort, the ring language $\LL_{\rm Rings}$ for
the residue field, and the Presburger language $\LL_{\rm PR}$ for
the value group together with maps $\ord$ and $\ac$ as in
\ref{dscf}. By taking the product of the measure $\mu^m$ with the
counting measure  on $k_K^n \times \ZZ^r$ one defines a measure
still denoted by $\mu^m$ on $K^m \times k_K^n \times \ZZ^r$.

One defines the dimension of a
$\cL'$-definable subset $X$ of $K^m\times k_K^n\times\ZZ^r$ as the
dimension of its projection  $\pi (X) \subset K^m$. If $X$ is of dimension $d$,
one defines a measure $\mu^d$ on $X$ extending the previous construction  on $X$ by setting
\begin{equation}
\mu^d (W) := \int_{\pi (X)}\pi_! (\11_W) \mu^d
\end{equation}
with $\pi_! (\11_W)$ the function $y \mapsto {\rm card}(\pi^{-1}(y)
\cap W)$.

For such an $X$, one defines $\cC (X)$ as the $\QQ$-algebra of
functions on $X$ generated by functions $\alpha$ and $p^{- \alpha}$
with $\alpha : X \rightarrow \ZZ$ definable in $\cL'$. Note that
this definition coincides with the previous one when $n = r = 0$.
Since $\cL'$ is interpretable in $\cL$, the formalism developed in
this section extends to $\cL'$-definable objects in natural way. In
particular the definitions of $\cC^{\leq d}$, $C^d$, ${\rm I}C$,
$\vert \Jac \vert$, $\delta$, $\cC^{\rm exp}$, etc, extend readily
to $\cL'$-definable objects and we have:

\begin{theorem}\label{alwaysthesame}\begin{enumerate}
\item[(1)]The statement of Theorem \ref{mainK} extends to $\Def (\cL')$
after adding the additional axiom:
\medskip

\noindent\textup{(A5)}\, Let $\pi : X \times
 k_K^n
  \times \ZZ^r
\rightarrow X$ be the projection with $X$ in $\Def (\cL')$. For any
$\varphi$ in ${\rm I}C (  X \times k_K^n \times \ZZ^r)$ and every
$x$ in $X$,
$$
\pi_! (\varphi) (x) =
\sum_{\pi (y) = x} \varphi (y).
$$
\item[(2)]
The statement of Theorem \ref{mainKexp} extends to $\Def (\cL')$.
\end{enumerate}
\end{theorem}

\subsubsection{}
Fix $\Lambda$ in $\Def (\cL')$. We consider the category
$\Def_{\Lambda}(\cL')$ whose objects are $\cL'$-definable morphisms
$f:S\to \Lambda$, a morphism $g : (f:S\to \Lambda) \rightarrow
(f':S'\to \Lambda)$ being a morphism $g : S \rightarrow S'$ in $\Def
(\cL')$ with $f' \circ g = f$. For $f:S\to \Lambda$ in
$\Def_{\Lambda}(\cL')$ we define the relative dimension of $S$ over
$\Lambda$ as the maximum of the dimension of the fibers of $f$. For
$d$ in $\ZZ$, define $\cC^{\leq d}(S\to\Lambda)$ as the ideal of
$\cC(S)$ generated by the characteristic functions of
$\cL'$-definable subsets of $S$ of relative dimension $\leq d$ over
$\Lambda$. Set
\begin{equation}C^d(S\to\Lambda):=\cC^{\leq
d}(S\to\Lambda)/\cC^{\leq d-1}(S\to\Lambda),
\end{equation}
 and define the graded group
\begin{equation}
C(S\to\Lambda):=\oplus_d C^d(S\to\Lambda).
\end{equation}
For every $\lambda$ in $\Lambda$ there exists a graded group
homomorphism called restriction to $\lambda$
\begin{equation}
|_{f^{-1}(\lambda)}:C(S\to\Lambda)\to C(f^{-1}(\lambda)),
\end{equation}
sending $\varphi$ in
$C(S\to\Lambda)$ to its restriction to the fiber $f^{-1}(\lambda)$.

We define ${\rm I} C(S\to \Lambda)$ as the graded subgroup of
$C(S\to\Lambda)$ consisting of $\varphi\in C(S\to\Lambda)$ such
that, for every $\lambda$ in $\Lambda$, the restriction $\varphi|_{f^{-1}(\lambda)}$ lies in
${\rm I} C(f^{-1}(\lambda))$.

One defines similarly
$\cC^{\leq
d}(S\to\Lambda)^{\rm exp}$,
$C(S\to\Lambda)^{\rm exp}$,
${\rm I} C(S\to \Lambda)^{\rm exp}$ and
\begin{equation}
|_{f^{-1}(\lambda)}:C(S\to\Lambda)^{\rm exp}\to C(f^{-1}(\lambda))^{\rm exp}.
\end{equation}

If $g : S \rightarrow S'$ is an isomorphism in $\Def_{\Lambda}
(\cL')$ between subsets of relative dimension $d$, on denotes by
$\vert \Jac_{\Lambda} g \vert$ the function in $C^d (S \rightarrow
\Lambda)$ such that
\begin{equation}
\vert \Jac_{\Lambda} g\vert_{|{f^{-1}(\lambda)}}
=
\vert \Jac (g_{|{f^{-1}(\lambda)}}) \vert
\end{equation}
for every $\lambda$ in $\Lambda$.

\begin{prop}\label{relKint}
For $g:S\to S'$ a morphism in  $\Def_{\Lambda} (\cL')$, there exists a unique
morphism
\begin{equation}\label{eqr}
g_{!\Lambda}:{\rm I} C(S\to\Lambda)\to {\rm I} C(S'\to\Lambda)
\end{equation}
which sends $\varphi\in {\rm I} C(S\to\Lambda)$ to the unique $\psi\in
 {\rm I} C(S'\to\Lambda)$ such that for each $\lambda\in\Lambda$
\begin{equation}
(g|_{S_{\lambda}})_!(\varphi|_{f^{-1}(\lambda)}) =
\psi|_{S'_{\lambda}},
\end{equation}
with $S_\lambda$ and $S'_\lambda$ the fibers and
$(g|_{S_{\lambda}})_!$ the direct image as constructed above,
and similarly a morphism
\begin{equation}\label{eqre}
g_{!\Lambda}:{\rm I} C(S\to\Lambda)^{\rm exp}\to {\rm I} C(S'\to\Lambda)^{\rm exp}.
\end{equation}
Furthermore, these morphisms $g_{!\Lambda}$ satisfy the relative
analogues of properties \textup{(A1)-(A4)} of Theorem \ref{mainK},
\textup{(A1)-(A4)} of Theorem \ref{mainKexp}, and (A5) of Theorem
\ref{alwaysthesame} respectively,
 where in \textup{(A4)} of Theorem
\ref{mainK}, $\vert \Jac\vert$ is replaced by its relative analogue
$\vert \Jac_{\Lambda} \vert$.
\end{prop}

\begin{proof}One can either remark that the proofs of
Theorem \ref{mainK}, Theorem \ref{mainKexp}
and Theorem \ref{alwaysthesame} carry over litterally
to the relative case, or deduce it from the absolute case
using
Proposition \ref{linear}
and Proposition \ref{linearexp}.
\end{proof}

Since ${\rm I}C(\Lambda \to\Lambda) = \cC (\Lambda)$, when $g$ is
the morphism $S\to \Lambda$, one gets from (\ref{eqr}) a morphism
\begin{equation}\label{eqrt}
\mu_{\Lambda} :{\rm I} C(S\to\Lambda)\to \cC (\Lambda).
\end{equation}

\section{Specialization and transfer}\label{transfert}

In this section we obtain new results on specialization to $p$-adic
and $\FF_q \llp t \rrp$-integration and a transfer principle for
exponential integrals with parameters from $\QQ_p$ and from $\FF_q
\llp t \rrp$. Some of the results which are announced in
\cite{miami} are generalized here to exponential constructible
functions. The specialization principle given here generalizes the
one of \cite{JAMS}.

\subsection{Specialization to valued  local  fields}\label{nd}

\subsubsection{Notation}Let $k$ be a number field with ring of integers $\cO$.We denote by
$\cA_{\cO}$ be the collection of all the $p$-adic completions of $k$
and all finite field extensions of $k$. We denote by $\cB_{\cO}$ the
set of all local fields of positive characteristic over $\cO$, that
is, endowed with an $\cO$-algebra structure. For $N>0$, denote by
$\cC_{\cO,N}$ the collection of $K$ in $\cA_{\cO}\cup \cB_{\cO}$
with ${\rm char}K>N$ and write $\cC_{\cO}$ for $\cC_{\cO,1}$. By
$\cA_{\cO,N}$, resp.~$\cB_{\cO,N}$, denote $\cC_{\cO,N}\cap
\cA_{\cO}$, resp.~$\cC_{\cO,N}\cap \cB_{\cO}$.

For $K$ in $\cC_{\cO}$, we write $R_K$ for its valuation ring, $M_K$
for the maximal ideal, $k_K$ for its residue field, and $q(K)$ for
the number of elements of $k_K$. For each choice of a uniformizer
$\varpi_K$ of $R_K$, there is a unique multiplicative map $\ac:
K^\times \to k_K^\times$ which extends the projection $R_K^\times\to
k_K^\times$ and sends $\varpi_K$ to $1$, and we extend this by
seting $\ac(0)=0$.
We denote by $\cD_K$ the collection of additive characters
$\psi:K\to\CC^\times$ such that
$$
\psi(x)=\exp((2\pi i / p) {\rm Tr}_{k_K}(\bar x))
$$ for $x\in R_K$, with $p$ the characteristic of $k_K$, ${\rm Tr}_{k_K}$
the trace of $k_K$ over its prime subfield and $\bar x$ the natural
projection modulo $M_K$ of $x$ into $k_K$.

\subsubsection{Interpretation of functions}
As a language that can be intepreted in all the fields of $\cC_\cO$,
we shall use $\LO:=\LPre (\cO \llb t\rrb )$, that is, the language
$\LPre$ with coefficients in $k$ for the residue field sort and
coefficients in $\cO \llb t\rrb $ for the valued field sort.
(Instead of $\cO \llb t\rrb$, any subring of $\cO \llb t\rrb$
containing $\cO[t]$ can be used as coefficients.)
 To say that a definable subassignment is defined in
the language $\LO$, we say that it belongs to $\Def(\LO)$, and for a
constructible function we say likewise that it belongs to $\cC (S,
\LO)$, $\cC (S, \LO)^{\rm exp}$, and so on, when it is defined in
$\LO$.

For every uniformizer $\varpi_K$ of $R_K$, one may consider $K$ as an
$\cO\llb t\rrb $-algebra via the morphism
\begin{equation}
\lambda_{\cO,K}:\cO\llb t\rrb \to K:\sum_{i\in\NN}a_it^i\mapsto
\sum_{i\in\NN}a_i\varpi_K^i.
\end{equation}
Hence, if one interprets elements $a$ of $\cO\llb t\rrb $ as
$\lambda_{\cO,K}(a)$, an $\cO\llb t\rrb $-formula $\varphi$ defines
for all $K$ in $\cC_\cO$ a definable subset $\varphi_K$ of $K^m \times
k_K^n \times \ZZ^r$ for some $m$, $n$, $r$, for every choice of uniformizer
$\varpi_K$ of $R_K$.

On the other hand, the formula $\varphi$ gives rise to a definable
subassignment $X$ of $h[m,n,r]$ and if $\varphi'$ gives rise to the
same subassignment $X$ then $\varphi_K=\varphi'_{K}$ for all $K$ in
$\cC_{\cO,N}$ for some large enough $N$, independently of the choice
of uniformizer.\footnote{This follows from either from Ax and Kochen
\cite{AK1}, \cite{AK3}, Er{\v s}ov \cite{Ersov}, Cohen \cite{cohen},
Pas \cite{Pas}, or others, or from a small variant of Proposition
5.2.1 of \cite{JAMS} (a result of Ax-Kochen-Er{\v s}ov type that
uses ultraproducts and follows from the Theorem of Denef-Pas).}

With a slight abuse of notation, for $X$ a definable subassignment
of $h [m, n, r]$ in $\Def(\LO)$, we write $X_K$ for the subassigment defined by
$\varphi_K$
where $\varphi$ is a $\Def(\LO)$-formula defining $X$, which is
well determined for $K$ in $\cC_{\cO,N}$ for some large enough $N$,
as explained above. Similarly, if $f:X\to Y$ is a $\LO$-definable
morphism, we obtain a function $f_K:X_K\to Y_K$ for all $K$ in
$\cC_{\cO,N}$ for some $N$.

With a similar abuse of notation, we can interpret a function $\varphi$ in
$\cC(X,\LO)$ as a function $X_{K}\to \QQ$, for $N$ large enough and
$K$ in $\cC_{\cO,N}$, as follows.

First suppose that $\varphi$ is in $K_0(\RDef_X (\LO))$ and of the
form $[\pi:W\to X]$ for some $\LO$-definable subassignment $W$ in
$\RDef_X(\LO)$. For $K$ in $\cC_\cO$, consider $W_{K}$, which is a
subset of $X_{K} \times (k_K)^{\ell}$ for some $\ell$, and consider
the natural projection $\pi_K:W_K\to X_K$. Then one sets
\begin{equation}\varphi_K:
\begin{cases}X_K\to\QQ& \\
x\mapsto {\rm card} \left(\pi_K^{-1}(x)\right).&
\end{cases}
\end{equation}
Similarly as before, this is well determined for $N$ large enough
and $K$ in $\cC_{\cO,N}$. By linearity that construction extends to $K_0(\RDef_X (\LO))$.

Let us now define $\varphi_K$ when $\varphi$ lies in $\cP(X)$. If one expresses $\varphi$
in terms of  $\LL$ and of definable morphisms $\alpha:X\to\ZZ$,
replacing $\LL$ by
$q_K$ and  each $\alpha$ by
$\alpha_K:X_K\to\ZZ$, one gets a function $\varphi_K:X_K\to\QQ$
again well determined for $K$ in $\cC_{\cO,N}$ when  $N$ is large enough.
By tensor product, this defines $\varphi_K$ for general
$\varphi$ in $\cC(X,\LO)$.

Next we interpret $\varphi$ in $\cC(X,\LO)^{\rm exp}$ as a function
$\varphi_{K,\psi_K}:X_{K}\to \CC$,
for $K$ in $\cC_{\cO,N}$ when  $N$ is large enough and for every  $\psi_K$ in $\cD_K$, as follows.

First suppose that $\varphi$ in $K_0(\RDef_X (\LO))^{\rm exp}$ is
of the form $[W,g,\xi]$ with $W$ a $\LO$-definable subassignment,
where $g:W\to h[1,0,0]$ and $\xi:W\to h[0,1,0]$  are $\LO$-definable
morphisms. For $K$ in $\cC_\cO$, consider $W_{K}$, $g_K:W_K\to K$,
and $\xi_K:W_K\to k_K$, and consider the projection $\pi:W_K\to
X_K$. Then, for $\psi_K$ in $\cD_K$, one sets
\begin{equation}\varphi_{K,\psi_K}:
\begin{cases}X_K\to\QQ & \\
 x
\mapsto \sum_{y\in \pi_K^{-1}(x)}\psi_K(g_K(y))\exp((2\pi i /p){\rm
Tr_{k_K}}(\xi_K(y))). &
\end{cases}
\end{equation}

Similarly as before, this is well determined for $N$ large enough
and for all $K$ in $\cC_{\cO,N}$ and all $\psi_K$ in $\cD_K$.
The construction being compatible  with the previous one, this
defines $\varphi_{K,\psi_K}$ for general $\varphi\in
\cC(X,\LO)^{\rm exp}$, by tensor product.

\subsubsection{Integration}
Let $K$ be in $\cC_\cO$ and consider a - not necessarily definable - subset  $A$ of $K^m\times
k_K^n \times \ZZ^r$. Let $A'$ be the
image of $A$ under the projection $K^m\times k_K^n \times \ZZ^r\to
K^m$ and define the dimension of $A$ as the dimension of the Zariski
closure of $A'$ in $\AA_K^m$ with  $\dim \emptyset:=-1$. Let $f : A
\rightarrow \Lambda$ be any function, with $\Lambda$ a subset of
$K^{m'} \times k_K^{n'} \times \ZZ^{r'}$.  The relative dimension
of $f$ is defined to be the maximum of the dimensions of the
fibers; it is also called the relative dimension of $A$ over
$\Lambda$.

For $A$ and $A'$ as above, we denote by  $\bar A$ be the Zariski closure of
$A'$ in $\AA_K^m$. If $\bar A$ is of dimension $d$, we consider
the canonical $d$-dimensional measure $\mu^d$ on $\bar A (K)$, cf.
\cite{serre}, \cite{oesterle}, and put the counting measure on $k_K^n \times
\ZZ^r$. We shall still write $\mu_d$
for the product measure on $\bar A (K)\times k_K^n \times \ZZ^r$,
formed by taking the product of the above $\mu_d$ with the
counting measure, and we finally still denote by $\mu^d$ its restriction to $A$,
similarly as in section \ref{sec:padics}.

We denote by  $\cF (A)$ the 
algebra  of all functions $A
\rightarrow \CC$. 
and we say a
function $\varphi$ in $\cF (A)$ is integrable in dimension $d$ if
$A$ and $\varphi$ are measurable and $\varphi$ is integrable with
respect to the measure $\mu^d$. More generally, we say that a
function $\varphi\in \cF (A)$ is integrable in dimension $e$ if
the support $B$ of $\varphi$ is of dimension $e$ and the
restriction $\varphi|_{B}$ is integrable in dimension $e$ as
defined above.
For $e \geq 0$ an integer, we denote by $\cF^{\leq e}(A)$ the ideal
of $\cF (A)$ of functions with support of dimension $\leq e$ and we
set
\begin{equation}F^e (A):=
\cF^{\leq e}(A) / \cF^{\leq e - 1}(A)\quad
 \text{and}
 \quad F (A) := \oplus_e F^e (A).\end{equation}
We define ${\rm I} F^e (A)$ as the subgroup of $F^e (A)$ consisting
of functions in $F^e (A)$ which are integrable in dimension $e$ and denote by
$\mu : {\rm I} F^e (A) \rightarrow \CC$ as the integration
operator. We set
${\rm I} F (A):=\oplus_e {\rm I} F^e (A)$
and extend $\mu$ to
$\mu : {\rm I} F (A) \rightarrow \CC$
by linearity.

Let $f : A \rightarrow \Lambda$ be a mapping as before. Let
$\cF^{\leq e} (A\rightarrow \Lambda)$ be the ideal of $\cF (A)$ of
functions with support of relative dimension $\leq e$ over
$\Lambda$. The groups $F^e (A\rightarrow \Lambda)$ and $F
(A\rightarrow \Lambda)$ are defined correspondingly.
 For every  $\lambda$ in $\Lambda$, there is a natural restriction map, which is a graded group
 homomorphism,
\begin{equation}
|_{f^{-1} (\lambda)}:F_K (A \rightarrow \Lambda)\to F_K
(f^{-1}(\lambda)),
\end{equation}
defined by sending $\varphi$ in $F_K (A \rightarrow \Lambda)$ to the
restriction of $\varphi$ to the fiber $f^{-1}(\lambda)$.
We define ${\rm I} F_K (A \rightarrow \Lambda)$ as the graded
subgroup of $F_K (A \rightarrow \Lambda)$ of Functions whose
restrictions to all fibers lie in ${\rm I} F$, where restriction
is as just defined, and we  denote by $\mu_{\Lambda}$ the unique
mapping
\begin{equation}\mu_{\Lambda}:{\rm I} F (A \rightarrow \Lambda) \rightarrow \cF
(\Lambda)
\end{equation}
such that $ \mu_{ \Lambda} (\varphi) (\lambda) = \mu
(\varphi|_{f^{-1} (\lambda)})$ for every $\varphi$ in ${\rm I} F
(A \rightarrow \Lambda)$ and every point $\lambda$ in $\Lambda$.

\subsubsection{}
We still have to go one step further in the interpretation of functions. Let
$f : S \rightarrow \Lambda$ be a morphism in $\Def (\LO)$. Let
$\varphi$ be in $C (S \rightarrow \Lambda, \LO)^{\rm exp}$. The
Function $\varphi$ is the class of a tuple $(\varphi_d)_d$ with
$\varphi_d$ in $\cC^{\leq d} (S \rightarrow \Lambda, \LO)^{\rm
exp}$, where only finitely many components are nonzero. Then, for
$N>0$ large enough, $K\in\cC_{\cO,N}$, and $\psi_K\in \cD_K$, each
function $\varphi_{d,K,\psi_K}$ lies in $\cF^{\leq d}_{K}
(S_{K})$, and by taking the class of $(\varphi_{d,K,\psi_K})_d$,
we get a function $\varphi_{K,\psi_K}$ in $F_{K} (S_{K}
\rightarrow \Lambda_{K})$.

The following result says that the motivic exponential integral
specializes to the corresponding integrals over the local fields
of high enough residue field characteristic.

\begin{theorem}[Specialization Principle]\label{compres3} Let
$f : S \rightarrow \Lambda$ be a morphism in $\Def (\LO)$. Take
$\varphi$ in ${\rm I}C (S \rightarrow \Lambda, \LO)^{\rm exp}$. Then there exists
$N>0$ such that for all $K$ in $\cC_{\cO,N}$, every choice of a
uniformizer $\varpi_K$ of $R_K$, and all
$\psi_K$ in $ \cD_K$, 
$\varphi_{K,\psi_K}$ lies in ${\rm I}F_{K} (S_{K}
\rightarrow \Lambda_{K})$ and 
\begin{equation*}
\left(\mu_\Lambda(\varphi)\right)_{K,\psi_K} = \mu_{\Lambda_K}(\varphi_{K,\psi_K}).
\end{equation*}
\end{theorem}
\begin{proof}
Let us first consider the case  where $\varphi$ lies in ${\rm I}C (S \rightarrow \Lambda, \LO)$.
We can assume
$\varphi$ lies  in ${\rm I}C_+ (S \rightarrow \Lambda, \LO)$, using  notations from
\cite{cl}. In \cite{cl}, the definition of relative integrability of $\varphi$ and the value of the relative integral  were defined simultaneously
along the following lines.
One may assume $S$ is a definable subassigment of $\Lambda [m, n, r] :=\Lambda \times  h [m, n, r]$
and, using cell decomposition and induction, it is enough, by Theorem 14.1.1 of \cite{cl}
to consider the behaviour
of the integrability condition and the computation of the integral for: 1) projection along $\ZZ$-variables,
2) projection along residue field variables, 3) projections
$\Lambda [m, n, r]  \rightarrow \Lambda [m- 1, n, r] $ when $S$ is a $0$-cell,
4) projections
$\Lambda [m, n, r]  \rightarrow \Lambda [m- 1, n, r] $ when $S$ is a $1$-cell adapted to $\varphi$.
Note that certainly given $\varphi$ the cell decompositions involved here will specialize
to cell decomposition defined by the specialized conditions when $N$ is large enough.
This is a special instance of the compactness argument in model theory.
In 1), one can assume $\varphi$ is a   Presburger function, that is lies in $\cP_+ (S)$
with the notation of loc. cit. In that case, the integrability condition was built from the start to be compatible with speciailization, since it was expressed by ``summability when $\LL$ replaced by $q > 1$''. Also
the relative integral was defined by summing up series in powers of $\LL$ and specializes to summing over  $\ZZ^r$ with respect to
the counting measure. Step 2) is tautogically compatible with specialization.
In step 3) a function $\LL^{- \ordjac f }$, defined almost everywhere occurs,
and for $N$ large enough it specializes to $\vert \Jac f_K \vert$. By the change of variables formula for integrals over
fields in $\cC_{\cO}$, it follows that 3) is compatible with specialization.
Finally step 4) is compatible with specialization since the relative motivic volume of a $1$-cell $Z$
specializes
to the volume of the corresponding $Z_K$, for $N$ large enough, by definition.

When $\varphi$ lies in $\cC (S,\LO)^{\rm exp}$, the statement about
compatibility of relative integrability with specialization holds by
the previous construction. The construction of the relative integral
of $\varphi$ can be performed along similar lines as before.
Specialization for steps 1), 2) and 3) holds for the same reasons as
before and only step
4)
needs to be considered. It follows from our constructions than it
is enough to show that the relative integral of the function $E (z)$
on a $1$-cell with special coordinate $z$ specializes to the
corresponding one over $K$, for $N$ large enough, which is clear by
construction, cf.
also
Lemma \ref{hypertriv}.
\end{proof}

\subsection*{}
For $p$-adic fields, we can say more, using the formalism of
\ref{sr}.

\begin{theorem}[Specialization Principle]\label{compres0}
Let $\Lambda$ be in $\Def (\LO)$ and let $f : S \rightarrow S'$ be a
morphism in $\Def_{\Lambda} (\LO)$. Let $\varphi$ be in ${\rm I}C (S
\rightarrow \Lambda, \LO)^{\rm exp}$. Then there exists $N>0$ such
that for all $K$ in $\cA_{\cO,N}$, each choice of a uniformizer
$\varpi_K$ of $R_K$, and all $\psi_K$ in $\cD_K$ the Function
$\varphi_{K,\psi_K}$ lies in ${\rm I}C (S_{K}\rightarrow
\Lambda_{K})^{\rm exp}$ and such that
\begin{equation*}
\left(f_{!\Lambda}(\varphi)\right)_{K,\psi_K}
=
f_{!\Lambda_K}(\varphi_{K,\psi_K}).
\end{equation*}
\end{theorem}
\begin{proof}
Similar to the proof of Theorem \ref{compres3}.
\end{proof}

\subsection{Transfer principle for integrals with parameters}

We start by proving the following abstract form of the Transfer principle:

\begin{prop}\label{strong} Let $\varphi$ be in
$\cC (\Lambda, \LO)^{\rm exp}$. Then, there exists an integer $N$
such that for all $K_1,K_2$ in $\cC_{\cO,N}$ with $k_{K_1}\simeq
k_{K_2}$ the following holds:
\begin{gather*}
\varphi_{K_1,\psi_{K_1}} = 0 \mbox{ for all }
\psi_{K_1}\in\cD_{K_1}
\\
\text{if and only if}\\
 \varphi_{K_2,\psi_{K_2}} = 0\mbox{ for all  }
\psi_{K_2}\in\cD_{K_2}.
\end{gather*}
\end{prop}

\begin{proof}
We first consider the case when  $\varphi$ lies in $\cC(\Lambda, \LO)^{\rm e}$.
Suppose that $\Lambda$ is a $\LO$-definable subassignment of
$h[m,n,r]$. We give a proof by induction on $m$. For $m=0$, the proof goes as follows.
By quantifier elimination, any finite set of formulas needed
to describe $\varphi$ can be taken to be valued field quantifier
free. It follows that
\begin{equation}\varphi_{K_1}=\varphi_{K_2}\end{equation}
for $K_1$ and $K_2$ in $\cC_{\cO,N}$ with $k_{K_1}\simeq k_{K_2}$
and $N$ large enough, since two ultraproducts $K = \prod_{\cU} K_i$
and $K '= \prod_{\cU} K'_i$ of fields  $K_i$ and $K'_i$ in $\cC_\cO$
with $k_{K_i}\simeq k_{K_i'}$ over a nonprincipal ultrafilter $\cU$
on a set $I$ are elementary equivalent, as soon as $K$ and $K'$ have
characteristic zero.

Now assume $m>0$. By 
applying inductively
the Cell Decomposition Theorem \ref{np}, we can construct a
$\LO$-definable morphism
\begin{equation}
f:\Lambda\to
h[0,n',r']
\end{equation}
for some $n'$, $r'$, and $\tilde \varphi\in \cC(h[0,n',r'],
\LO)^{\rm e}$, such that $\varphi=f^*(\tilde \varphi)$. Necessarily,
$\tilde \varphi$ is unique. By the induction hypothesis,
\begin{equation}\tilde \varphi_{K_1}=0
\text{ if and only if } \tilde\varphi_{K_2}=0\end{equation}
for $K_1$ and $K_2$ in $\cC_{\cO,N}$ with $k_{K_1}\simeq k_{K_2}$
and $N$ large enough.  Since $\varphi_K=f_K^*(\tilde \varphi_K)$ for $K$ in
$\cC_{\cO,N}$ when $N$ is large enough,
the result follows for general $m$ and for $\varphi$ in
$\cC (\Lambda, \LO)^{\rm e}$.

In general, when $\varphi$ lies in $\cC
(\Lambda, \LO)^{\rm exp}$, we write
$\varphi$ as a finite sum of the form
\begin{equation}\label{tcl}
\sum_{i=1}^\ell E(g_i)e(\xi_i)[X_i\to\Lambda]\varphi_i,
\end{equation}
with $\varphi_i\in \cC (\Lambda, \LO)$.

After  finitely partitioning $\Lambda$, we may suppose that there is a
partition of $\{1,\ldots,\ell\}$ into parts $B_r$ such that
 \begin{equation}\ord (g_i(x_i)-g_j(x_j))<0
 \end{equation} for
all $i\in B_{r_1}$, all $j\in B_{r_2}$, all $r_1\not=r_2$, all
$\lambda\in \Lambda$ and all $x_i\in X_i$, $x_j\in X_j$ lying above
$\lambda$, and such that \begin{equation} \ord (g_i(x_i)-g_j(x_j))\geq 0\end{equation} for all
$i,j\in B_{r}$, all $r$, all $\lambda\in \Lambda$ and all $x_i\in
X_i$, $x_j\in X_j$ lying above $\lambda$.

\begin{claim}There exists $N>0$ such that for all
$K$ in $\cC_{\cO,N}$ the statement
\begin{equation}
\varphi_{K,\psi_K}=0\mbox{ for every }\psi_K\in\cD_K
\end{equation}
is equivalent to
\begin{equation}\label{tclb}
\sum_{i\in B_r}
(e(\xi_i)[X_i\to\Lambda]\varphi_i)_{K,\psi_K}=0\mbox{ for every $r$
and for every }\psi_K\in\cD_K.
\end{equation}
\end{claim}

Since the left hand side of (\ref{tclb}) is in fact independent of
the choice of character $\psi_K$, the proposition directly follows
from the claim and the treatment of the case $\varphi$ in
$\cC(\Lambda, \LO)^{\rm e}$.

Let us now prove the claim. By compactness there exists $N_0$ such that
for all $K\in\cC_{\cO,N_0}$ the partition $\{B_r\}_r$ satisfies
the following property: for
every $i\in B_{r_1}$, every $j\in B_{r_2}$, every $r_1\not=r_2$, every
$\lambda\in \Lambda_K$ and every  $x_i\in X_{iK}$, $x_j\in X_{jK}$
lying above $\lambda$
 \begin{equation}\ord (g_{iK}(x_i)-g_{jK}(x_j))<0,
 \end{equation} and,
 for every  $i,j\in B_{r}$, every
$r$, every  $\lambda\in \Lambda_K$ and every $x_i\in X_{iK}$, $x_j\in
X_{jK}$ lying above $\lambda$,
 \begin{equation}
 \ord (g_{iK}(x_i)-g_{jK}(x_j))\geq 0.
 \end{equation}
Now the claim follows from Lemma \ref{triveq}.
\end{proof}

\begin{lem}\label{triveq} Let $K$ be in $\cC_{\cO}$.
Let $c_i$ be in $\CC$ and $x_i\in K$ with $\ord(x_i-x_j)< 0$ for
$i\not=j$, $i,j=1,\ldots,n$. For every $\psi$ in $\cD_K$ consider the
exponential sum
\begin{equation*}
S_\psi:=\sum_{i=1}^n c_i\psi(x_i).
\end{equation*}
Suppose that $S_\psi=0$ for all $\psi$ in $\cD_K$. Then $c_i=0$
for all $i$.
\end{lem}

\begin{proof}
We shall perform an induction on
$m:= -\min_i (\ord(x_i))$. If $m=0$
there is nothing to prove. So let us assume $m\geq 1$.

For every $n \geq 0$, we denote by $\cD_K (n)$ the set of
restrictions of the characters in $\cD_K$ to the ball
$\varpi_K^{-n}R$. We denote by $p$ the characteristic of $k_K$ and
we set $\pi := \varpi_K$ if $K$ is of characteristic $p$ and $\pi :=
p$ if $K$ is of characteristic $0$. We fix elements $y_1$, \dots,
$y_r$ of $\varpi_K^{-m} R$ whose images in $\varpi_K^{-m} R / \pi
\varpi_K^{-m} R$ form a $\FF_p$-basis. For $a = (a_1, \dots, a_r)$
in $\{0, \dots, p - 1\}^r$, we denote by $B_a$ the ball $\ord (x -
\sum_j a_j y_j) \geq \ord(\pi \varpi_K^{-m})$. Let us fix $\psi_0$
in $\cD_K (m - 1)$. There are exactly $p^r$ characters in $\cD_K
(m)$ extending $\psi_0$. Indeed, such characters are determined by
their value on $y_1$, \dots, $y_r$, hence if we denote by $\zeta_{j,
i}$, for $1 \leq i \leq p$, the $p$ distinct complex numbers such
that $\zeta_{j, i}^p = \psi_0 (p y_j)$, they are  in one to one
correspondence with the set of tuples $(\zeta_{j, i})_j$, via $\psi
\mapsto (\psi (y_j))$.

We may rewrite
$S_{\psi}$
as
\begin{equation}
S_{\psi}
= \sum_{a \in \{0, \dots, p - 1\}^r}
\prod_{1 \leq j \leq r} \psi (y_j)^{a_j} S_{a , \psi_0}
\end{equation}
with
\begin{equation}
S_{a , \psi_0}
=
\sum_{x_i \in B_a} c_i \psi_0 \Bigl(x_i -   \sum_{1 \leq j \leq r}  a_j y_j\Bigr).
\end{equation}

For fixed $j$, the $p \times p$-matrix
$A_j := (\zeta^{\ell}_{j, i})_{i, \ell}$, $0 \leq \ell \leq p - 1$, is an invertible Vandermonde matrix.
It follows that the
Kronecker (tensor) product
matrix $A_1 \otimes \dots \otimes A_r$
with coefficients
$\prod_{1 \leq j \leq r} \zeta_{j, i}^{\ell_j}$, $0 \leq \ell_j \leq p - 1$,
is an
invertible
$p^r \times p^r$-matrix.
Thus, the vanishing of $S_{\psi}$ for every $\psi$ in
$\cD_K (m)$ implies the vanishing of all the sums
$S_{a , \psi_0}$ for every $\psi_0$ in
$\cD_K (m - 1)$, and the induction hypothesis allows to conclude.
\end{proof}

Now we can prove the following fundamental
transfer principle for exponential integrals:

\begin{theorem}[Transfer principle for exponential integrals]\label{strongaxkc}
Let $S \rightarrow \Lambda$ and $S' \rightarrow \Lambda$ be
morphisms in $\Def (\LO)$. Let $\varphi$ be in ${\rm I}C (S
\rightarrow \Lambda, \LO)^{\rm exp}$ and $\varphi'$ in ${\rm I}C
(S' \rightarrow \Lambda, \LO)^{\rm exp}$. Then, there exists an
integer $N$ such that for all $K_1,K_2$ in $\cC_{\cO,N}$ with
$k_{K_1}\simeq k_{K_2}$ the following holds
\begin{gather*}
\mu_{ \Lambda_{K_1}} (\varphi_{K_1,\psi_{K_1}}) = \mu_{
\Lambda_{K_1}} (\varphi'_{K_1,\psi_{K_1}}) \mbox{ for all }
\psi_{K_1}\in\cD_{K_1}
 \\
 \text{if and only if}
 \\
\mu_{ \Lambda_{K_2}} (\varphi_{K_2,\psi_{K_2}}) = \mu_{
\Lambda_{K_2}} (\varphi'_{K_2,\psi_{K_2}}) \mbox{ for all  }
\psi_{K_2}\in\cD_{K_2}.
\end{gather*}
\end{theorem}

\begin{proof}
By taking the
disjoint union of $S$ and $S'$ over $\Lambda$ and linearity,
it is enough to prove the following particular case of the result:
if $S \rightarrow \Lambda$ is a morphism in $\Def (\LO)$ and
$\varphi$ is in ${\rm I}C (S \rightarrow \Lambda, \LO)^{\rm exp}$,
 there exists an integer $N$ such that for all $K_1,K_2$ in
$\cC_{\cO,N}$ with $k_{K_1}\simeq k_{K_2}$ the following holds:
\begin{gather*}
\mu_{ \Lambda_{K_1}} (\varphi_{K_1,\psi_{K_1}}) = 0 \mbox{ for all
}
\psi_{K_1}\in\cD_{K_1}\\
\text{if and only if}\\
 \mu_{ \Lambda_{K_2}} (\varphi_{K_2,\psi_{K_2}}) = 0 \mbox{ for all  }
\psi_{K_2}\in\cD_{K_2},
\end{gather*}
which follows directly from Theorem \ref{compres3} and Proposition
\ref{strong}.
\end{proof}

\begin{remark}\label{rl}
Without exponentials, a form of Theorem \ref{strongaxkc} can already
be found in \cite{miami}. As mentioned in the introduction, it
should have a wide range of applications to $p$-adic representation
Theory and the Langlands Program, in particular to various forms of
the Fundamental Lemma. Let us note in our approach  there is no need
any more to make assumptions of local constancy, as was done in
\cite{CH}. Let us recall that the Fundamental Lemma for unitary
groups has been proved by  Laumon and Ng\^o
 \cite{ln} over functions fields and that  Waldspurger deduced the case of $p$-adic fields \cite{wal}
 by representation theoretic
techniques. A typical situation when Theorem
\ref{strongaxkc} can be directly applied  is the Jacquet-Ye conjecture \cite{JY}, a relative
version of the Fundamental Lemma involving integrals of additive characters, which has been proved
by
Ng\^o \cite{Ngo} over
functions fields and by Jacquet \cite{J} in general.
\end{remark}

\subsection*{Acknowledgment}
{\small During the realization of this project, the first author was
a postdoctoral fellow of the Fund for Scientific Research - Flanders
(Belgium) (F.W.O.) and was supported by The European Commission -
Marie Curie European Individual Fellowship with contract number HPMF
CT 2005-007121.}

\bibliographystyle{amsplain}

\end{document}